\documentclass[preprint,12pt]{elsarticle}

\usepackage{graphics}
\usepackage{graphicx}
\usepackage{epsfig}
\usepackage{xcolor,soul}
\usepackage{anysize}
\usepackage{array}
\usepackage[fleqn]{amsmath}
\usepackage{latexsym,amssymb,amscd,amsthm,amsxtra}
\usepackage{subfigure}
\usepackage{tabularx}
\usepackage{enumerate}
\usepackage{enumitem}
\usepackage{multirow}
\usepackage{float}
\usepackage{placeins}
\usepackage{booktabs}
\usepackage{times}
\usepackage{soul,xcolor}
\usepackage{algorithm}
\usepackage{algpseudocode}
\usepackage[permil]{overpic}
\usepackage{tikz}
\usetikzlibrary{shapes.geometric, arrows}
\usetikzlibrary{arrows.meta}
\usetikzlibrary{positioning}

\tikzstyle{startstop} = [rectangle, rounded corners, minimum width=2.5cm, minimum height=0.8cm, text centered, draw=black, fill=red!30, font=\small]
\tikzstyle{process} = [rectangle, minimum width=3.5cm, minimum height=0.8cm, text centered, draw=black, fill=blue!30, font=\small, text width=3.5cm, align=center]
\tikzstyle{arrow} = [thick,->,>=stealth]

\newtheorem{remark}{Remark}[section]

\newtheorem{lemma}{Lemma}[section]

\usepackage[most]{tcolorbox}

\usepackage{arydshln}

\begin{document}
\begin{frontmatter}

 \title{B-spline-Based ALE-MFS Framework for Evolving Domains}

 \author{Muhammad Ammad\corref{cor1}\fnref{label1}}
\ead{21481199@life.hkbu.edu.hk}
\cortext[cor1]{Corresponding author}
\author{Leevan Ling\fnref{label1}}
\author{Shu Ma\fnref{label1}}
%\cortext[cor2]{}

\affiliation[label1]{organization={Department of Mathematics, Hong Kong Baptist University},
             city={Kowloon Tong, Kowloon},
             country={Hong Kong}}

 \begin{abstract}
We develop and analyze a B-spline based arbitrary Lagrangian-Eulerian method of fundamental solutions (ALE-MFS) for curvature-driven motion of two-dimensional evolving domains. Boundary points move with the material to track the geometric flow, while interior points move within an ALE framework via a harmonic extension of the boundary velocity, computed by a meshless MFS with sources on a fixed auxiliary circle, thus avoiding volumetric meshing. Boundary normals and curvature are reconstructed by an adaptive local B-spline scheme that remains robust for strongly nonconvex shapes and large deformations. A posteriori error estimates are obtained from a hatmatrix formulation of leave-one-out cross-validation (LOOCV) for both square collocation and zero-padded least-squares systems, and are complemented by maximum principle indicators for harmonic problems. Numerical experiments on circular, star-shaped, and amoeba-like domains show that square collocation suffices for moderately complex geometries, while zero-padded least-squares significantly improves interior velocity regularity and pointwise transport accuracy for strongly nonconvex shapes, without altering the source or collocation sets. The ALE-MFS algorithm also generates high-quality moving meshes for ALE-finite element methods, with larger minimum angles and slower mesh ratio growth than classical FEM mesh-motion strategies, suggesting a practical and easily integrable alternative for challenging moving-interface simulations.
% \bigskip
% \noindent \textbf{Research highlights}:
% \begin{itemize}
%     \item  Localized tensor-product B-spline patches enable flexible, efficient interpolation on evolving point clouds without global meshing.
%     \item Interpolant/control points are updated directly at each time step, avoiding reinterpolation.
%     \item Control-point coefficients have a clear geometric meaning, allowing in-place updates under curvature-driven motion.
%     \item The conditioning-aware formulation for local interpolation preserves the quality of the interpolation under surface changes.
%     \item Adaptive point management (insertion/removal) maintains appropriate sampling and surface resolution.
%     \item The numerical results confirm a reliable simulation of curvature-driven flows and coupled surface–scalar interactions.
% \end{itemize}

 \end{abstract}


\end{frontmatter}

\section{Introduction}\label{sec:intro}
\begin{figure}
    \centering
    \includegraphics[width=0.3\textwidth]{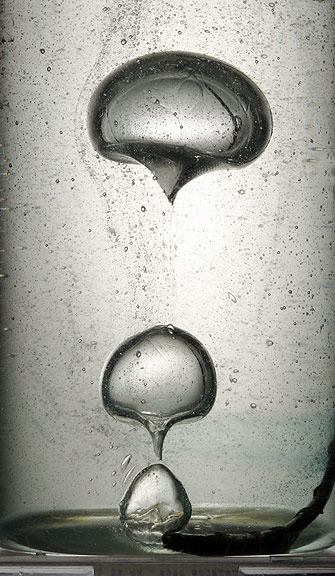} 
    \caption{Illustration of air bubbles rising and deforming in a liquid column, a classic moving boundary problem.}
    \label{fig:moving_boundary}
\end{figure}

Accurate numerical simulation of moving boundary problems remains a central challenge in applied mathematics and computational physics, with applications ranging from fluid dynamics to biomechanics and materials science \cite{guus2023+moving,guo2021+moving,verzicco2023immersed}. Such problems are characterized by domains that evolve over time due to underlying physical processes (see Figure~\ref{fig:moving_boundary} for an example of rising air bubbles), which require computational methods that can efficiently and accurately track these dynamic boundaries. Traditional fixed-mesh methods often struggle to capture evolving interfaces, leading to significant interest in more flexible numerical frameworks \cite{garcia2025+MFF,masud2007+meshmove}.

Among these frameworks, the arbitrary Lagrangian-Eulerian (ALE) approach has emerged as a powerful paradigm for simulating moving boundary problems \cite{subramaniam2013ale,kovacs2019ale}. The ALE method bridges the gap between the purely Lagrangian approach, where mesh distortion becomes problematic for large deformations \cite{pila2019introduction}, and the Eulerian approach that requires additional interface-capturing techniques \cite{mehta2022level}. By decoupling the motion of the mesh from the flow of the material, ALE maintains the accuracy of the interface tracking while preserving the quality of the mesh through controlled node displacement \cite{takizawa2022mesh,farinatti2004mesh}.

In practice, ALE is performed most frequently with finite elements (ALE-FEM) \cite{kovacs2018female,lan2021female,edelmann2022femALe}. However, ALE-FEM faces significant difficulties in complex geometries. First, the computational cost escalates rapidly because FEM requires discretization of the entire domain; this burden becomes especially acute in three dimensions. Frequent mesh regeneration and solution interpolation further increase overhead, particularly when extreme deformations or topological changes occur \cite{li2025optimal,gaburro2021unified}. Second, maintaining mesh quality over long simulations is challenging: large deformations can degrade element quality, and dynamically adding or removing points to adapt to evolving geometries often necessitates intricate remeshing algorithms and re-computation. Although adaptive remeshing strategies \cite{vautrin2020automatic} and machine learning-enhanced mesh motion \cite{haubner2024learning} have alleviated some issues, fundamental bottlenecks remain in mesh-based implementations.

These limitations motivate our investigation of meshless alternatives, particularly the method of fundamental solutions (MFS). As a boundary-type technique, MFS eliminates volumetric meshing by representing solutions as aggregates of fundamental solutions \cite{cheng2025introduction}. Unlike FEM, MFS only requires collocation points on the boundary, significantly reducing the number of degrees of freedom and computational expense, especially for large domains. The method has demonstrated remarkable accuracy in boundary value problems \cite{cheng2020overview}, while its meshless nature provides inherent advantages for evolving domains by avoiding the need for mesh updates or topological adjustments during domain evolution \cite{chantasiriwan2009+131}. Recent extensions have addressed historical limitations by coupling with particular solutions \cite{chen2012+144} and time-stepping schemes \cite{valtchev2008time}. 

To overcome these limitations, we develop a hybrid ALE–MFS framework that retains ALE’s interface-tracking capabilities while using the meshless structure of MFS. The approach eliminates repeated, domain-wide remeshing and associated interpolation costs, and mitigates mesh distortion by relying on fundamental-solution bases, thus enabling accurate simulations in complex, time-evolving geometries \cite{alves2004+17}. The early results highlight clear benefits for complex interfacial dynamics \cite{karageorghis2023efficient} and free-surface flows \cite{xiao2017+932}.

The remainder of this paper is structured as follows. Section~\ref{sec:model} introduces the ALE framework for moving boundary problems, establishing the mathematical and geometric formulations that underpin our approach. Section~\ref{sec:mfs} details the method of fundamental solutions, including its adaptation to domains with prescribed or inferred boundary velocities, with particular attention to the integration of curvature-driven flows via B-spline interpolation. In Section~\ref{sec:error}, we discuss \emph{a posteriori} error indicators suited to square and overdetermined MFS systems, providing the necessary tools for rigorous error assessment. Section~\ref{sec:nu_ex_5} presents an extensive suite of numerical experiments in a range of geometries and configurations, systematically analyzing the interplay between discretisation, error quantification, and geometric complexity. Building on these experiments, Section~\ref{sec:ale_fem_comparison} generates comparable moving meshes based on the proposed B-spline ALE--MFS algorithm and compares the mesh quality with the classical ALE--FEM. Finally, Section~\ref{sec:conclusion} summarizes the main findings and discusses the prospects for future developments in meshless ALE methodologies for moving interface problems.

\section{Model Formulation via the ALE Framework}\label{sec:model}
\begin{figure}[ht!]
    \centering
    \begin{tikzpicture}
        \node[anchor=south west, inner sep=0] (img) at (0,0) {\includegraphics[width=0.2\textwidth]{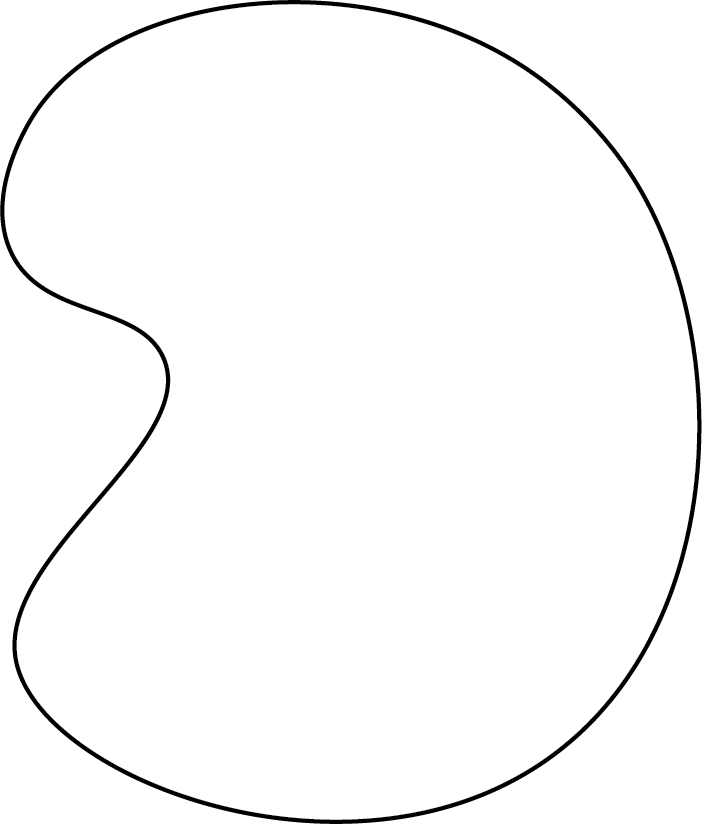}};
        \begin{scope}[x={(img.south east)}, y={(img.north west)}]
            \node at (0.55,0.5) {$\Omega(t)$};
            \node at (-0.1,0.88) {$\Gamma(t)$};
            \draw[->, thick, black] (0.95,0.7) -- (1.1,0.8) node[right] {$\vec{V}_\Gamma$};
        \end{scope}
    \end{tikzpicture}
    \caption{Annotated moving domain \(\Omega(t)\) and boundary \(\Gamma(t)\) in \(\mathbb{R}^2\) with boundary velocity \(\vec{V}_\Gamma\).}
    \label{fig:ale-domain-annotated}
\end{figure}
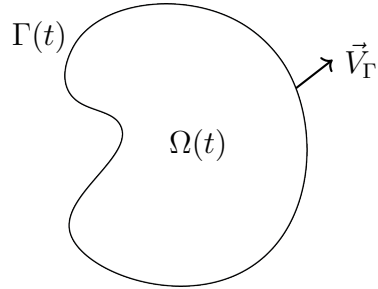

In the context of moving boundary problems, the ALE framework provides a robust approach to handling evolving domains by decoupling the motion of the computational mesh from the material flow. Let \(\Omega(t) \subset \mathbb{R}^d\) denote the spatial domain dependent on time \(t \geq 0\), with its limit defined as \(\Gamma(t) = \partial \Omega(t)\). For illustrative purposes, a two-dimensional representation (\(d=2\)) of the domain and its boundary is provided in Figure~\ref{fig:ale-domain-annotated}.
The evolution of the domain is described through a family of bijective mappings, referred to as the ALE mappings, given by
\[
\mathcal{A}_t: \overline{\Omega}_0 \to \overline{\Omega}(t),
\]
where \(\Omega_0 = \Omega(0)\) represents the initial (reference) configuration of the domain at \(t=0\). The trajectory of a point in the domain is thus expressed as
\[
\vec{x} = \mathcal{A}_t(\vec{X}), \quad \vec{X} \in \Omega_0,
\]
where \(\vec{X}\) denotes the coordinates in the reference configuration and \(\vec{x}\) represents the corresponding coordinates in the current configuration at time \(t\). The associated mesh velocity field, which governs the motion of the computational mesh, is defined as
\begin{equation}
\vec{u}(\vec{x}, t) := \left. \frac{\partial}{\partial t} \mathcal{A}_t(\vec{X}) \right|_{\vec{X} = \mathcal{A}_t^{-1}(\vec{x})}
\label{eq:ALEframework}
\end{equation}

In moving boundary problems, the normal velocity of the boundary, denoted by \(\vec{V}_\Gamma\), is often prescribed based on the underlying physical process (e.g., fluid-structure interaction or surface evolution). This boundary velocity is imposed on \(\Gamma(t)\); however, extending this velocity into the interior of \(\Omega(t)\) to define a mesh velocity field \(\vec{u}\) is not unique and requires a specific strategy. Various approaches exist in the literature, including algebraic prolongation \cite{allen2002algebraic}, elastic deformation models \cite{aygun2023physics}, and harmonic extensions \cite{WuCai2014}, each offering different trade-offs in terms of computational cost and mesh quality preservation.

For this study, we adopt the harmonic extension due to its ability to produce a smooth velocity field that minimizes mesh distortion. Specifically, the velocity of the mesh \(\vec{u}\) is determined as the solution to the following boundary value problem:
\begin{equation}
\begin{cases}
\Delta \vec{u} = 0 & \text{in } \Omega(t), \\
\vec{u} = \vec{V}_\Gamma & \text{on } \Gamma(t),
\end{cases}
\label{eq:harmonic-extension}
\end{equation}
where $\Delta$ represents the component-wise Laplacian operator. 
Since $\vec{u} = \vec{V}_\Gamma$ on $\Gamma(t)$, $\vec{u}$ and $\vec{V}_\Gamma$ generate the same evolution of the boundary $\Gamma(t)$. The harmonic extension ensures that the prescribed boundary velocity is exactly matched on $\Gamma(t)$, while the interior mesh velocity is the smoothest possible extension in the $H^{1}$-sense. This property is particularly advantageous for maintaining the quality of the mesh over extended simulations, making harmonic extension a widely adopted choice in ALE-based computations for moving-boundary problems.

Although the mesh velocity and its harmonic extension are central to traditional mesh-based ALE implementations, our ultimate goal is to transition to a meshless framework, specifically the MFS, in subsequent sections. The discussion here establishes the geometric and analytical foundation of the moving domain and velocity field, which will inform the adaptation of ALE principles into a meshless context.

\section{Formulation of the Method of Fundamental Solutions}\label{sec:mfs}
The MFS provides a meshless computational framework for solving partial differential equations (PDEs), offering significant advantages for moving boundary problems due to its ability to handle evolving geometries without the burden of volumetric meshing. Within the scope of our study, we apply MFS to the Laplace equation in a time-dependent domain \(\Omega(t) \subset \mathbb{R}^2\), with a boundary \(\Gamma(t) = \partial \Omega(t)\), as defined in Section~\ref{sec:model}. This section outlines the MFS formulation tailored to the Laplace equation in the context of moving boundaries under the arbitrary Lagrangian-Eulerian (ALE) framework.

\subsection{MFS with Prescribed Boundary Velocities}
\label{subsec:mfs_prescribed}

\begin{figure}[pt!]
    \centering
    \begin{tikzpicture}
        \node[anchor=south west, inner sep=0] (img) at (0,0) {\includegraphics[width=0.3\textwidth]{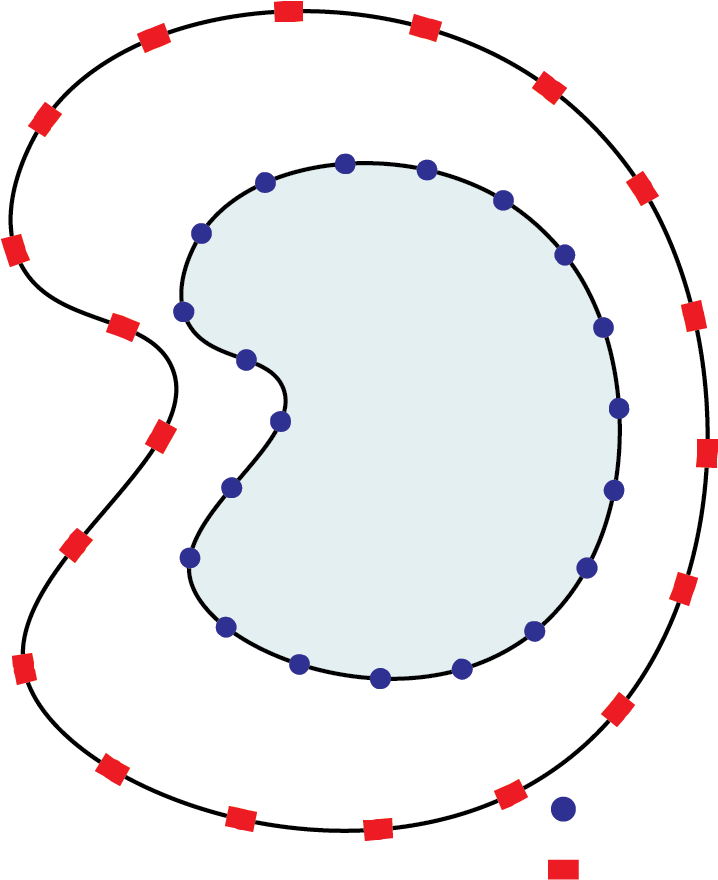}};
        \begin{scope}[x={(img.south east)}, y={(img.north west)}]
            \node at (0.55,0.5) {$\Omega(t)$};
             \node at (1.2,0.075) {Collocation points};
              \node at (1.12,0.005) {Source points};
            %\node at (-0.1,0.88) {$\Gamma(t)$};
            \draw[->, thick, black] (0.86,0.5) -- (1.03,0.59) node[right] {Main boundary};
             \draw[->, thick, black] (0.92,0.74) -- (1.,0.78) node[right] {Pseudo boundary};
        \end{scope}
    \end{tikzpicture}
    \caption[Schematic of the MFS setup]{Schematic of the MFS setup: the domain $\Omega$ (shaded), boundary collocation points $N_c$ (blue dots), and auxiliary source points $N_s$ (red square) located outside the domain.}
    \label{fig:mfs_schematic}
\end{figure}
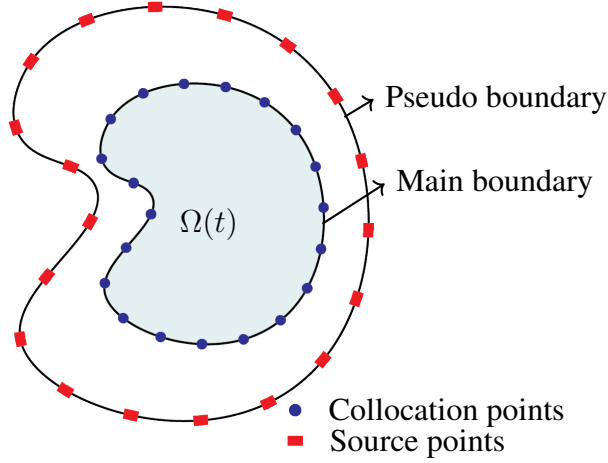

Suppose that the boundary velocity is prescribed by the underlying physics as
\[
g(\vec{x}, t) = \vec{V}_\Gamma(\vec{x}, t), \quad \vec{x} \in \Gamma(t).
\]
We seek a vector-valued field \(\vec{u}(\vec{x}, t)\) that satisfies the Laplace equation with Dirichlet boundary conditions.
\begin{equation}
\label{eq:mfs-continuous}
\Delta \vec{u}(\vec{x}, t) = \vec{0} \quad \text{in } \Omega(t), \qquad \vec{u}(\vec{x}, t) = g(\vec{x}, t) \quad \text{on } \Gamma(t).
\end{equation}
Unlike traditional methods, no volumetric mesh is introduced; instead, MFS circumvents the need for volumetric meshing by expressing the solution as a superposition of fundamental solutions, each singular at an auxiliary source point located outside the domain (see Figure~\ref{fig:mfs_schematic}).
Explicitly, at a fixed time \(t\), the MFS ansatz for each component of the vector-valued field \(\vec{u}(\vec{x}, t)\) takes the form of a
\begin{equation}
\label{eq:mfs-ansatz}
u^{(k)}(\vec{x}, t) = \sum_{j=1}^{N_s} \alpha_j^{(k)}(t) \, \Phi(\vec{x}, \vec{y}_j), \quad k=1,2,
\end{equation}
where \(\{\vec{y}_j\}_{j=1}^{N_s}\) are fixed source points outside \(\overline{\Omega(t)}\), \(\{\alpha_j^{(k)}(t)\}\) are unknown coefficients to be determined for each component \(k\), and \(\Phi(\vec{x}, \vec{y}_j)\) denotes the fundamental solution of the Laplacian in \(\mathbb{R}^2\):
\begin{equation}
\label{eq:fundamental-solution}
\Phi(\vec{x}, \vec{y}_j) = -\frac{1}{2\pi} \log |\vec{x} - \vec{y}_j|.
\end{equation}
The enforcement of the boundary condition at a set of collocation points \(\{\vec{x}_i\}_{i=1}^{N_c} \subset \Gamma(t)\) produces a linear system for the coefficients \(\{\alpha_j^{(k)}(t)\}\), thus defining the MFS approximation. Specifically, imposing the Dirichlet condition from \eqref{eq:mfs-continuous} on each \(\vec{x}_k\) results in
\[
\sum_{j=1}^{N_s} \alpha_j^{(k)}(t) \, \Phi(\vec{x}_k, \vec{y}_j) = g_k(\vec{x}_k, t), \quad k=1,2, \quad \text{collocation indices}\; i=1, \dots, N_c.
\]
In practice, one often sets \(N_s = N_c\), leading to a square linear system that can be solved for the coefficients. Once determined, the approximate mesh-velocity field \(\vec{u}^{N_s}(\vec{x}, t)\) can be evaluated anywhere in \(\Omega(t)\). This meshless approach adapts naturally to the evolving geometry of \(\Gamma(t)\) driven by \(\vec{V}_\Gamma\), avoiding the computational overhead of mesh regeneration inherent in finite element methods within the ALE framework.

\subsection{MFS with Inferred (Unprescribed) Boundary Velocities}
\label{subsec:mfs_unprescribed}

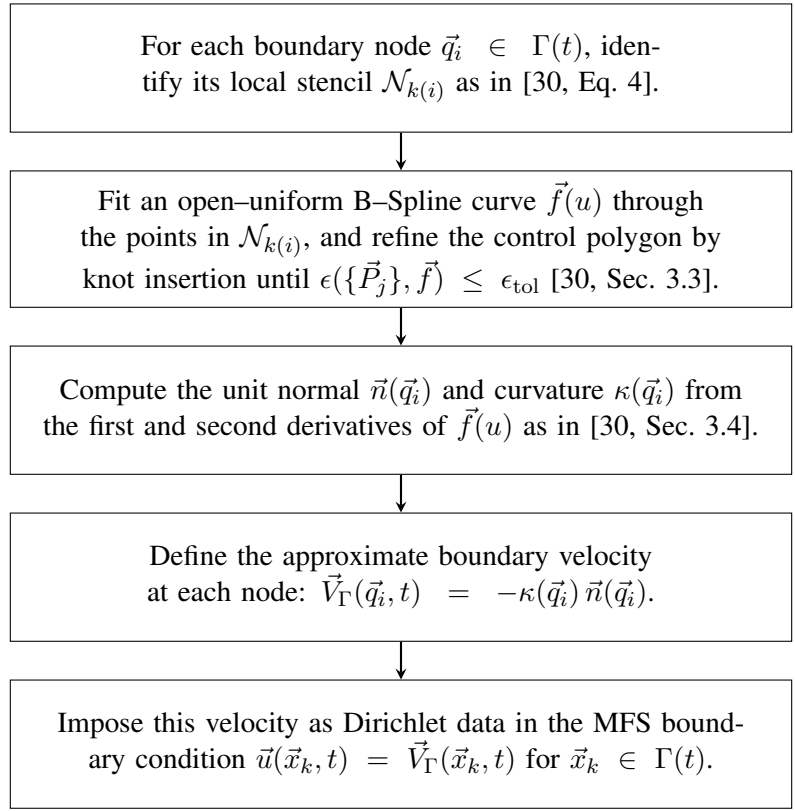
\begin{figure}[htbp]
\centering
\begin{tikzpicture}[
  node distance=1cm,
  box/.style={
    rectangle,
    draw=black,
    text width=10cm,
    minimum height=1.7cm,
    text centered,
    font=\small
  },
  arrow/.style={thick,->,>=stealth}
]

% Nodes
\node (step1) [box]
  {For each boundary node $\vec q_i\in\Gamma(t)$, identify its local stencil $\mathcal N_{k(i)}$ as in \cite[Eq.~4]{ammad2025eabe}.};

\node (step2) [box, below=0.5cm of step1]
  {Fit an open–uniform B–Spline curve $\vec f(u)$ through the points in $\mathcal N_{k(i)}$, and refine the control polygon by knot insertion until $\epsilon(\{\vec P_j\},\vec f)\le\epsilon_{\rm tol}$ \cite[Sec.~3.3]{ammad2025eabe}.};

\node (step3) [box, below=0.5cm of step2]
  {Compute the unit normal $\vec n(\vec q_i)$ and curvature $\kappa(\vec q_i)$ from the first and second derivatives of $\vec f(u)$ as in \cite[Sec.~3.4]{ammad2025eabe}.};

\node (step4) [box, below=0.5cm of step3]
  {Define the approximate boundary velocity at each node: $\vec V_\Gamma(\vec q_i, t) = -\kappa(\vec q_i)\,\vec n(\vec q_i)$.};

\node (step5) [box, below=0.5cm of step4]
  {Impose this velocity as Dirichlet data in the MFS boundary condition $\vec u(\vec x_k, t) = \vec V_\Gamma(\vec x_k, t)$ for $\vec x_k\in\Gamma(t)$.};

% Arrows
\draw [arrow] (step1) -- (step2);
\draw [arrow] (step2) -- (step3);
\draw [arrow] (step3) -- (step4);
\draw [arrow] (step4) -- (step5);

\end{tikzpicture}
\caption{Flowchart of the adaptive B–Spline procedure used to infer boundary velocities $\vec V_\Gamma$ from a discrete interface \cite[Sec.~3]{ammad2025eabe}.}
\label{fig:mfs_velocity_procedure}
\end{figure}

In many moving boundary problems, the normal velocity of the boundary,
\(\vec V_\Gamma\), is not directly prescribed by the physical model.
Instead, one must infer \(\vec V_\Gamma\) from geometric information,
typically a point cloud sampled from the interface. In this work, for
simplicity, we adopt the mean–curvature motion law
\(\vec V_\Gamma = -\,\kappa\,\vec n\) whenever an explicit physical
velocity is unavailable. To compute \(\kappa\) and \(\vec n\) from the
discrete boundary, we employ the adaptive B–Spline machinery of
\cite[Sec.~3]{ammad2025eabe}. The overall procedure is summarized
schematically in Fig.~\ref{fig:mfs_velocity_procedure}.

Using the local smoothness and endpoint–interpolation properties of
open–uniform B–Splines, this approach yields high–order accurate
approximations of \(\kappa\), \(\vec n\), and hence \(\vec V_\Gamma\),
without requiring any external velocity prescription. The resulting
boundary data are then seamlessly incorporated into the MFS system
\eqref{eq:mfs-continuous}, enabling a fully meshless, curvature–driven
evolution of the domain.

%%%%%%%%%%%%%%%%%%%%%%%%%%%%%%%%%%%%%%%%%%%%%%%

\subsection{Selection of Collocation and Source Points}\label{sec:points_mfs}
The choice of collocation and source sets is critical in the MFS.  It depends on both boundary conditions (harmonic and non-harmonic) and geometry (near-circular boundaries, strongly varying curvature, corners). Several practical deployments have been assessed in the literature; see, e.g., \cite{li2013mfs_nonharmonic,chen2016choosing}. In this study, we adopt a workflow suited to evolving domains and point cloud inputs.

When a boundary parametrization is available, we initialize collocation points from the parametrization. For general evolving interfaces (after the first time step), we work directly with boundary point clouds; local B-spline reconstruction \cite[Sec. 3]{ammad2025eabe} provides normals and curvature and allows resampling when needed. The method does not rely on globally uniform spacing.

In all experiments, we employ a classical, circle-based source set: sources are placed on a fixed circle strictly outside the domain. Concretely, if $\vec c$ denotes a chosen center (e.g., the centroid at $t=0$) and $R_\mathrm{s}>0$ is a radius larger than the maximum domain radius during the simulation, we set
\[
\vec y_j \;=\; \vec c \;+\; R_\mathrm{s}\,(\cos\theta_j,\;\sin\theta_j), 
\qquad \theta_j = \frac{2\pi(j-1)}{N_s}, \; j=1,\dots,N_s,
\]
and keep this circular source layer fixed as the domain evolves. This “sources-on-a-circle” placement is simple, robust, and avoids recomputing source normals or distances in moving geometries, which is advantageous in an ALE setting with point-cloud boundaries.

Alternative source placements such as boundary-tied dilation layers (homothetic) and normal-offset layers $\vec y_j=\vec x_j+d\,\vec n(\vec x_j)$ are well documented and can be effective, particularly on static or smoothly parameterized geometries \cite{chen2016choosing,li2013mfs_nonharmonic}. Because our domains evolve and we prioritize a stable, geometry-agnostic setup without repeated geometric preprocessing, we restrict attention to the fixed circular source layer in this study.

\section{Error Indicators}\label{sec:error}
Before presenting the numerical results, we outline the error indicators a~posteriori used to assess the accuracy of solutions obtained using the arbitrary Lagrangian--Eulerian method of fundamental solutions (ALE--MFS) for the Laplace equation in a time-dependent domain $\Omega(t) \subset \mathbb{R}^2$ with boundary $\Gamma(t)$.

Depending on the stabilization strategy, the boundary-collocation system for the MFS is treated in one of two square formulations :

\begin{itemize}
    \item \textbf{Classical square system}: The number of boundary collocation points is equal to the number of source points, which yields a square linear system. This is the standard MFS interpolation setting.
    \item \textbf{Zero-padded square system}: The same square system is augmented with zero rows to stabilize a least-squares solve while keeping the source and collocation sets unchanged. This construction preserves $N_c = N_s$ in the collocation matrix, but uses a least-squares perspective through zero padding. 
\end{itemize}

In this work, we use two complementary \emph{a~posteriori} indicators that apply to both the classical and zero-padded square settings: (i) the $\mathrm{LOOCV}\!-\!\ell_\infty$ indicator. For well-conditioned square interpolation systems, we exploit the classical fast LOOCV specialization. For general square systems, including the zero-padded least-squares formulation, we use the hat-matrix LOOCV formula, which remains valid and provides a rigorous pointwise \emph{a~posteriori} error indicator; (ii) the maximum principle boundary indicator, in which we monitor the maximal boundary error via the maximum principle, both as an external check and to contextualize the LOOCV bound.

\subsection{LOOCV Indicator}
To assess the predictive accuracy of the MFS approximation and to detect possible overfitting or ill-conditioning in the collocation system, we employ the Leave-One-Out Cross-Validation (LOOCV) error indicator, specifically in its $\ell_\infty$-norm variant. This a posteriori measure, rooted in statistical validation techniques, evaluates how well the model interpolates the boundary data by simulating the omission of each collocation point in turn. The concept of LOOCV traces back to early work in regression diagnostics and model selection, such as Allen's prediction sum of squares (PRESS) criterion \cite{allen1974relationship} and the generalized cross-validation by Golub et al. \cite{golub1979generalized}, with subsequent refinements in machine learning and approximation theory \cite{hastie2009elements}. In the context of meshfree methods such as RBF interpolation and MFS, it has been adapted as an efficient tool for parameter tuning and error estimation \cite{rippa1999algorithm,fasshauer2007meshfree}.

The general formulation relies on the hat matrix and applies uniformly to both square and overdetermined linear systems. Consider the collocation matrix $A \in \mathbb{R}^{N_c \times N_s}$ and the vector on the right-hand side $\vec{g}$, where the MFS coefficients are obtained as $\vec{\alpha} = A^{+} \vec{g}$, with $A^{+}$ denoting the Moore--Penrose pseudoinverse of $A$. The fitted values at the collocation points are then $\hat{\vec{g}} = A \vec{\alpha} = H \vec{g}$, where $H = A A^{+}$ is the hat matrix.

The LOOCV residual at the $j$-th collocation point is defined as
\[
e_j = \frac{g_j - \hat{g}_j}{1 - H_{jj}}, \qquad j = 1, \dots, N_c,
\]
which quantifies the error incurred when predicting $g_j$ using a model fitted to all other points. This expression can be efficiently computed for all $j$ without explicitly resolving the system $N_c$ times, using the properties of the hat matrix \cite{golub1979generalized,hastie2009elements}. For multi-component problems, such as those involving velocity fields with components $k=1,2$, we evaluate the residuals separately for each component and define the global LOOCV indicator as
\begin{equation}
E_{\mathrm{LOO},\infty} = \max_{k=1,2}\,\, \max_{1 \leq j \leq N_c} |e_j^{(k)}|. \label{eq:loocv_general}
\end{equation}

In the special case of a square, invertible system ($N_c = N_s$), the general formula simplifies to a computationally advantageous form introduced by Rippa \cite{rippa1999algorithm} for the RBF interpolation:
\begin{equation}
e_j^{(k)} = \frac{\alpha_j^{(k)}}{[A^{-1}]_{jj}}, 
\qquad j=1,\dots,N_s, 
\qquad \boldsymbol{\alpha}^{(k)} = A^{-1}\vec{g}^{(k)}.
\label{eq:rippa}
\end{equation}
This specialization takes advantage of the invertibility of $A$ to avoid the direct computation of the hat matrix, making it particularly efficient for well-conditioned square systems encountered in classical MFS settings.

For rectangular or singular systems (e.g., $N_c \neq N_s$), where the inverse $A^{-1}$ is unavailable, a heuristic algorithm, often termed the pseudoinverse--Rippa (pinv--Rippa) formula, extends Rippa's idea by replacing $A^{-1}$ with $A^{+}$:
\begin{equation}
\widetilde{e}_j^{(k)} = \frac{\alpha_j^{(k)}}{[A^{+}]_{jj}}, 
\qquad j=1,\dots,\min(N_s,N_c), 
\qquad \boldsymbol{\alpha}^{(k)} = A^{+}\vec{g}^{(k)}.
\label{eq:pinv-rippa}
\end{equation}
Although this coincides with \eqref{eq:rippa} in the square invertible case, it serves as an empirical approximation for overdetermined systems, producing a residual per source point (column of $A$) rather than per collocation point (row of $A$) \cite{rippa1999algorithm,fasshauer2007meshfree}. Unlike the hat-matrix LOOCV, which is rigorously derived, the pinv--Rippa algorithm lacks a formal theoretical justification in the rectangular setting and is used primarily as a practical diagnostic. We report its aggregated form
\begin{equation}
E_{\mathrm{pR},\infty} = \max_{k=1,2} \max_{1 \le j \le N_s} |\widetilde{e}_j^{(k)}|,
\label{eq:pinv-rippa-rect}
\end{equation}
alongside the general LOOCV indicator for comparative purposes.
%%%%%%%%%%%%%%%%%%%%%%%%%%%%%%%%%%%%%%%%%%%%%%%%%%%%%%%%%%%%%%%%%%%%%
\subsection{Maximum Principle Error Indicator}
To provide a direct measure of the solution error that is independent of the system configuration, we also employ the maximum principle for harmonic functions. This fundamental property asserts that the maximum and minimum values of a non-constant harmonic function on a bounded domain $\Omega$ are attained on the boundary $\Gamma$, thereby allowing boundary errors to control interior errors.

\begin{lemma}[Strong maximum principle; e.g., {\cite[Theorem~4, Section~2.2.3]{evans2022partial}}]
Let \(u \in C^2(\Omega)\cap C^0(\overline{\Omega})\) satisfy \(\Delta u = 0\) in a bounded domain \(\Omega \subset \mathbb{R}^2\). Then
\[
\max_{\overline{\Omega}} |u| = \max_{\partial \Omega} |u|.
\]
\end{lemma}

Let \(\{\vec{z}_p\}_{p=1}^{P} \subset \Gamma(t)\) be a set of boundary test points (typically with $P \geq N_c$ for dense sampling). For each component \(k=1,2\), define the pointwise error
\[
\epsilon_p^{(k)} = u_{N_s}^{(k)}(\vec{z}_p,t) - u_{\mathrm{true}}^{(k)}(\vec{z}_p,t), \qquad p=1,\dots,P,
\]
where $u_{N_s}^{(k)}$ denotes the MFS approximation with sources $N_s$, and $u_{\mathrm{true}}^{(k)}$ is the exact or reference solution. Motivated by the maximum principle, we employ the normalized boundary indicator
\begin{equation}
E_{\mathrm{MP},\infty} = \max_{k=1,2} \left(
\frac{\max_{1\le p \le P} |\epsilon_p^{(k)}|}
{\max_{1\le p \le P} |u_{\mathrm{true}}^{(k)}(\vec{z}_p,t)|}
\right), \qquad \vec{z}_p \in \Gamma(t),
\label{eq:mp-indicator}
\end{equation}
which provides a relative measure of the maximum boundary error and, by the lemma, bounds the interior error.

\section{Numerical Examples}\label{sec:nu_ex_5}
In this section, we conduct a comprehensive study of the numerical performance of the ALE-MFS on evolving planar domains. The boundary evolution is governed by the mean curvature flow, where the normal velocity at each boundary point is given by $\vec{V}_\Gamma = -\kappa \vec{n}$, with $\kappa$ denoting the curvature and $\vec{n}$ the unit outward normal vector. Both geometric quantities are computed at discrete boundary nodes via locally fitted B-spline interpolation, ensuring smooth and accurate representations even for complex and nonconvex geometries.

To investigate the robustness and accuracy of the ALE-MFS approach, we consider a sequence of representative examples exhibiting increasing geometric complexity. The examples are systematically organized into three parts, corresponding to distinct domain geometries and MFS system formulations, as summarized in Table~\ref{tab:numerical_examples_summary}:

\begin{table}[t!]
    \centering
    \renewcommand{\arraystretch}{1.2}
    \begin{tabular}{|c|c|c|c|} 
        \hline 
        \textbf{Part} & \textbf{Example} & \textbf{Geometry} & \textbf{MFS System} \\
        \hline\hline
        \multirow{1}{*}{Part I} 
            & 1 & Circle & Square  \\
        \hline\hline
        \multirow{3}{*}{Part II} 
            & 2 & Smooth Asterisk & Square  \\
        \cline{2-4}
            & 3 & Sharp-feature Asterisk & Square  \\
        \cline{2-4}
            & 4 & Amoeba & Square  \\
        \hline\hline
        \multirow{1}{*}{Part III}
            & 5 & Amoeba & Square (zero-padded) \\
            \hline
    \end{tabular}
    \caption{Summary of numerical experiments: geometry, and MFS system configuration used in each example.}
    \label{tab:numerical_examples_summary}
\end{table}
%& $E_{\mathrm{LOO},\infty}$ $(LOOCV-\ell_\infty)$
\begin{itemize}
    \item \textbf{Part I: Classical geometries.} This part is dedicated to canonical domains, such as the circle, serving as benchmarks for the ALE-MFS method in simple settings.
    \item \textbf{Part II: Non-smooth and complex geometries.} Here, we address domains with increasing geometric complexity and reduced regularity, including the smooth asterisk, sharp-feature asterisk, and amoeba shapes. These cases are designed to challenge the method's robustness in the presence of nontrivial boundary features.
\item \textbf{Part III: Overdetermined systems} In this part, we analyze the amoeba geometry using overdetermined MFS systems. We adopt the same source-point placement strategy used in the square-system experiments (circular source layer), without the additional normal-offset enhancement. Boundary accuracy is monitored using the LOOCV indicator, and we also report the coefficient-wise pinv--Rippa algorithm for rectangular systems.
\end{itemize}

For each case, numerical accuracy is rigorously assessed using the LOOCV indicator $E_{\mathrm{LOO},\infty}$, which provides a conservative upper bound on the pointwise boundary interpolation error. These diagnostics facilitate a systematic comparison of ALE-MFS performance across a variety of geometries and source placement strategies.

The results presented in the following subsections elucidate the interplay between domain geometry, source placement, and numerical accuracy, and highlight the strengths and limitations of ALE-MFS in the context of moving interface problems.

%%%%%%%%%%%%%%%%%%%%%%%%%%%%%%%%%%%%%%%%%%%%%%%%%%%%%%%%%%%%%%%%%%
\subsection{Comparison of MFS Parameters: Classical Circular Geometry}
In this subsection, we examine the sensitivity and accuracy of the method of fundamental solutions (MFS) under variation of discretization and parameter choices, using the circle as a canonical test case. The objective is to assess how fundamental parameters, such as the placement radius of source points and the boundary fill distance, influence the overall numerical stability and accuracy of the method.
\subsubsection*{Example 1: Circle Geometry with Uniformly Distributed Source Points}
We consider the initial domain as the unit circle $\Gamma(0)$, parameterized by $N$ boundary nodes. The MFS source points are uniformly placed in a concentric circle of radius $r_s > 2$, exterior to the physical boundary. The boundary evolves under mean curvature flow, governed by the velocity law up to a fixed time $t = T$.
\begin{figure}[t!]
\centering
\subfigure[Initial domain]{\includegraphics[width=3.in]{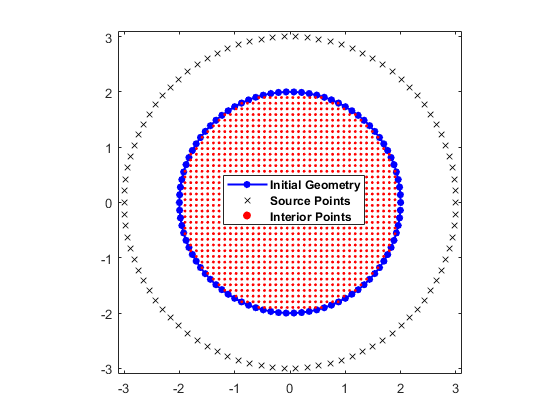}}\hspace{-1cm}
\subfigure[Evolved domain at $t = T$]{\includegraphics[width=3.in]{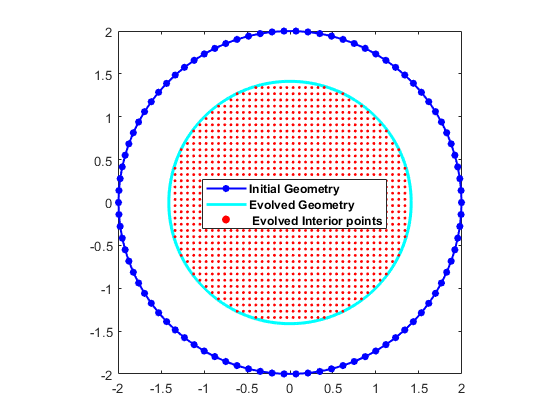}}\\
\subfigure[LOOCV and maximum principle error versus $r_s$]{
    \begin{overpic}[width=3in]{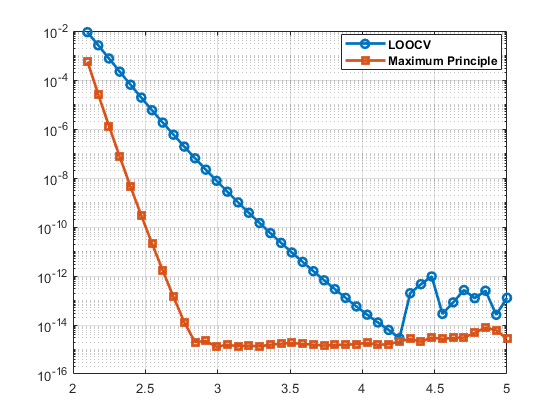}
        \put(20,240){\rotatebox{90}{\small $\ell_\infty$ norm error}}
        \put(500,10){$r_s$}
    \end{overpic}
}
\hspace{-1cm}
\subfigure[$E_{\mathrm{LOO},\infty}$ versus $r_s$ for various $h$]{
    \begin{overpic}[width=3in]{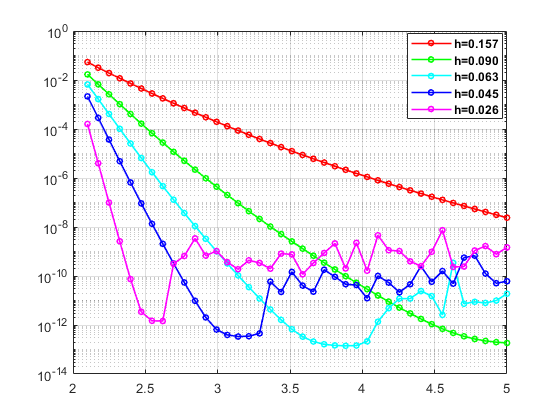}
        \put(20,220){\rotatebox{90}{\small $\ell_\infty$ norm error}}
        \put(500,10){$r_s$}
    \end{overpic}
}
\caption[Circular geometry with circular source points]{Circular geometry with circular source points.
(a) Initial geometry with uniform distribution of interior points and source points (on a concentric circle)
(b) Evolved boundary and interior points at \texttt{$t=T$}
(c) LOOCV indicator $E_{\mathrm{LOO},\infty}$ and maximum principle error $E_{\mathrm{MP}}$ versus $r_s$
(d) $E_{\mathrm{LOO},\infty}$ versus $r_s$ for various fill distances $h$.
}
\label{fig:example1}
\end{figure}
The primary aim of this experiment is to quantify the performance of the MFS as a function of the source circle radius $r_s$ and the boundary fill distance $h$ (the maximal chordal distance between consecutive boundary nodes). For error quantification, we focus on the LOOCV indicator $E_{\mathrm{LOO},\infty}$ and the maximum principle error $E_{\mathrm{MP}}$, both measured in the $\ell_\infty$ norm. As demonstrated in Figure~\ref{fig:example1}(c), $E_{\mathrm{LOO},\infty}$ consistently dominates $E_{\mathrm{MP}}$ across the entire range of source radii $r_s$, providing a sharp estimator for the boundary error in this canonical setting. The results further reveal that as the boundary fill distance $h$ is refined, $E_{\mathrm{LOO},\infty}$ exhibits uniform decay, indicative of the expected exponential convergence with respect to discretization density (see Figure~\ref{fig:example1}(d)). 

These observations confirm the reliability and practical effectiveness of $E_{\mathrm{LOO},\infty}$ as an error indicator, establishing a baseline for parameter selection and error assessment in more complex geometries considered in subsequent examples.

\subsection{Assessing LOOCV for MFS in Challenging Domains} \label{sec:part2}
After confirming the reliability of the LOOCV indicator for classical geometry, we proceed to examine its behavior in the context of non-smooth and highly non-convex domains. As the complexity of the boundary increases, several challenges emerge, including the degradation of interior point distribution and the potential breakdown of standard (square) MFS formulations. Through the following examples, we analyze the performance and limitations of the LOOCV error indicator in guiding parameter selection and assessing solution accuracy for increasingly challenging geometrical configurations.

\subsubsection*{Example 2: Smooth Asterisk Domain}
To advance our exploration beyond purely circular geometries, we consider as our initial domain $\Gamma(0)$ a smooth asterisk-shaped curve, a mild but non-trivial perturbation of the unit circle. The boundary is discretized with uniformly spaced nodes, and MFS source points are positioned on a concentric auxiliary circle with radius $r_s > 2$, as in the previous example. The evolution of the boundary is again governed by mean curvature flow, which, due to the smoothness of the initial shape, preserves a well-behaved interior point distribution throughout the simulation.

\begin{figure}[ht]
\subfigure[Initial domain]{\includegraphics[width=3in]{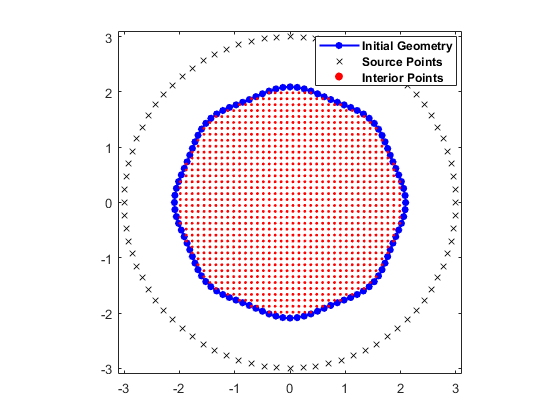}}
\subfigure[Evolved domain at $t = T$]{\includegraphics[width=3in]{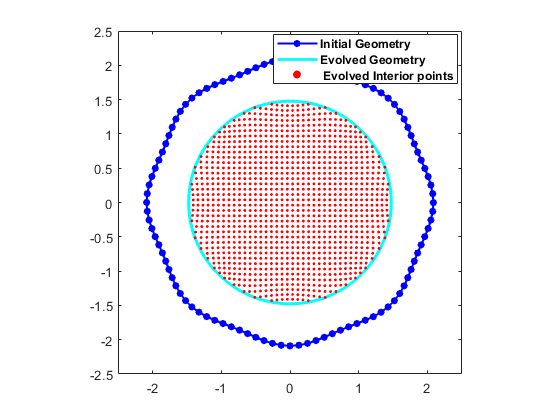}}\\
\subfigure[LOOCV and maximum principle error versus $r_s$]{
    \begin{overpic}[width=3in]{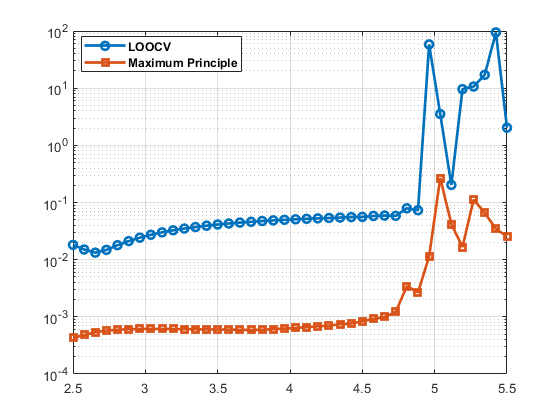}
        \put(20,240){\rotatebox{90}{\small $\ell_\infty$ norm error}}
        \put(500,10){$r_s$}
    \end{overpic}
}
\hspace{0.02in}
\subfigure[$E_{\mathrm{LOO},\infty}$ versus $r_s$ for various $h$]{
    \begin{overpic}[width=3in]{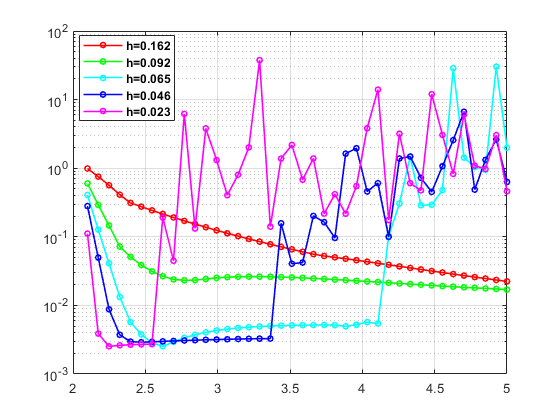}
        \put(20,220){\rotatebox{90}{\small $\ell_\infty$ norm error}}
        \put(500,10){$r_s$}
    \end{overpic}
}
\caption[Smooth asterisk geometry with circular source points]{Smooth asterisk geometry with circular source points
(a) Initial geometry with uniform distribution of interior points and source points (on a concentric circle)
(b) Evolved boundary and interior points at $t = T$
(c) LOOCV indicator $E_{\mathrm{LOO},\infty}$ and maximum principle error $E_{\mathrm{MP}}$ versus $r_s$
(d) $E_{\mathrm{LOO},\infty}$ versus $r_s$ for various fill distances $h$.
}
\label{fig:example2}
\end{figure}

In this setting, we systematically vary both the source circle radius $r_s$ and the boundary fill distance $h$, tracking their influence on the principal error indicators $E_{\mathrm{LOO},\infty}$ $(LOOCV-\ell_\infty)$ and $E_{\mathrm{MP}}$ (maximum principle error), both measured in the $\ell_\infty$ norm. The empirical findings in Figure~\ref{fig:example2}(c) demonstrate that $E_{\mathrm{LOO},\infty}$ continues to serve as a reliable overestimate of the maximum principle error across the explored parameter space. Notably, even as the domain departs from perfect circularity, $E_{\mathrm{LOO},\infty}$ retains its diagnostic sharpness, consistently bounding the true residual error.

Examining Figure~\ref{fig:example2}(d), we observe that the anticipated decay of $E_{\mathrm{LOO},\infty}$ with decreasing fill distance $h$ remains robust, confirming the stability of the MFS under moderate geometric perturbations. Although the convergence profile becomes somewhat less regular compared to the circle, the overall trend of exponential decay persists, and the error curves show only minor sensitivity to the increased geometric complexity.

These results illustrate that the LOOCV indicator not only adapts to mildly perturbed geometries but also continues to provide actionable guidance for parameter selection and error control. The successful extension of these behaviors to the smooth asterisk domain reinforces the versatility and predictive power of $E_{\mathrm{LOO},\infty}$ as an a posteriori error estimator.
\subsubsection*{Example 3: Sharp Feature Asterisk Domain}

In this example, the boundary $\Gamma(0)$ is initialized as an asterisk-shaped shape exhibiting moderately pronounced tips, features that are sharper than a smooth curve but do not reach the level of true geometric singularities. The MFS source points are distributed along a concentric exterior circle, and the domain is evolved according to mean curvature flow, which preserves the gently pointed structure throughout the simulation.

\begin{figure}[t!]
\subfigure[Initial domain]{\includegraphics[width=3in]{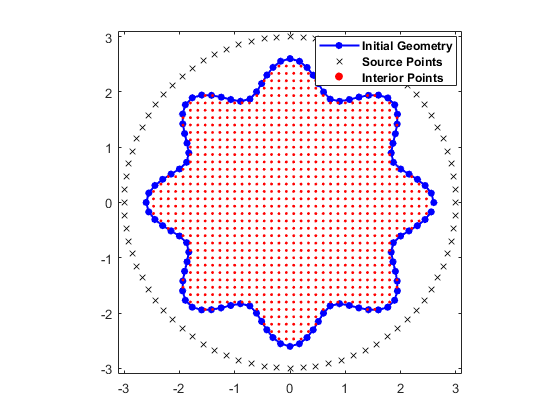}}
\subfigure[Evolved domain at $t = T$]
{  \begin{tikzpicture}
        \node[anchor=south west,inner sep=0] (image) at (0,0) {\includegraphics[width=3in]{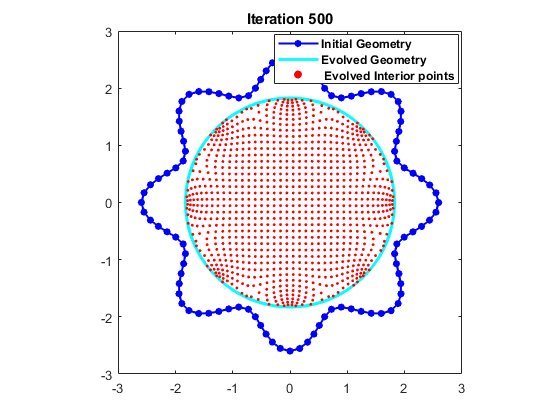}};
        \fill[white] (3.3,5.3) rectangle (4.9,5.8); 
      \end{tikzpicture}
    }\\
\subfigure[LOOCV and maximum principle error versus $r_s$]{
    \begin{overpic}[width=3in]{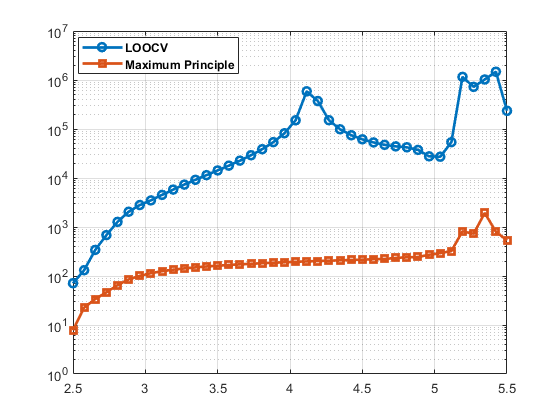}
        \put(20,240){\rotatebox{90}{\small $\ell_\infty$ norm error}}
        \put(500,10){$r_s$}
    \end{overpic}
}
\hspace{0.02in}
\subfigure[$E_{\mathrm{LOO},\infty}$ versus $r_s$ for various $h$]{
    \begin{overpic}[width=3in]{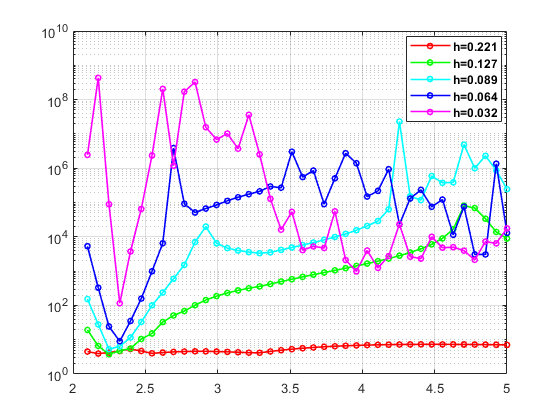}
        \put(20,220){\rotatebox{90}{\small $\ell_\infty$ norm error}}
        \put(500,10){$r_s$}
    \end{overpic}
}
\caption[Sharp feature asterisk geometry with circular source points]{Sharp feature asterisk geometry with circular source points
(a) Initial geometry with uniform distribution of interior points and source points (on a concentric circle)
(b) Evolved boundary and interior points at $t = T$
(c) LOOCV indicator $E_{\mathrm{LOO},\infty}$ and maximum principle error $E_{\mathrm{MP}}$ versus $r_s$
(d) $E_{\mathrm{LOO},\infty}$ versus $r_s$ for various fill distances $h$.
}
\label{fig:example3circle}
\end{figure}

These moderately acute tips introduce additional complexity compared to the previous examples. While not strictly corners in the mathematical sense, these features create regions of elevated curvature that challenge the numerical method. As illustrated in Figure~\ref{fig:example3circle}(b), the interior points remain evenly dispersed during evolution, avoiding problematic clustering or excessive sparsity.

Examining the error indicators, both $E_{\mathrm{LOO},\infty}$ and $E_{\mathrm{MP}}$ are found to be larger than in the smooth asterisk scenario, reflecting the increased difficulty of approximating boundary behavior in the vicinity of these tips. Nonetheless, the error values remain under reasonable control, and $E_{\mathrm{LOO},\infty}$ consistently bounds $E_{\mathrm{MP}}$ from above, confirming its reliability as an error gauge in this intermediate regime.

The results shown in Figure~\ref{fig:example3circle}(c,d) further reveal that, although the error curves exhibit more variability with respect to the source radius $r_s$ and fill distance $h$, the overall trend remains consistent with theoretical expectations. The heightened sensitivity is a natural consequence of the more intricate boundary geometry, yet the LOOCV indicator continues to provide valuable feedback for parameter selection and solution validation.

So, even when the domain features moderately sharp protrusions rather than true corners, the MFS paired with LOOCV maintains commendable robustness. This example demonstrates the method's capacity to handle boundaries with nontrivial curvature variations while still offering meaningful and practical error assessment.

\subsubsection*{Example 4: Amoeba Domain}
In this example, the initial domain is given by a highly non-convex, irregular “amoeba” shape, with MFS source points distributed along an exterior circle. As the boundary evolves, the interior point distribution deteriorates significantly, leading to pronounced clustering and sparsity, especially in regions of high curvature.

\begin{figure}[t!]
\subfigure[Initial domain]{\includegraphics[width=3in]{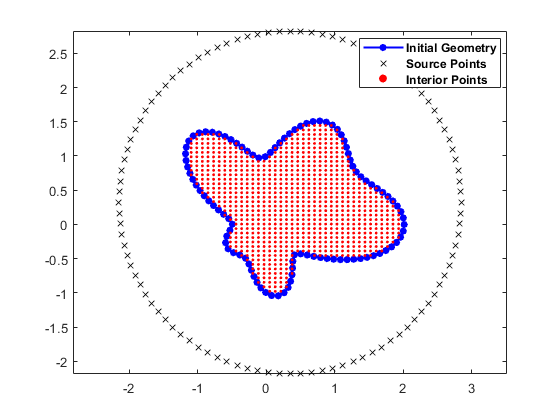}}
\subfigure[Evolved domain at few time step]{\includegraphics[width=3in]{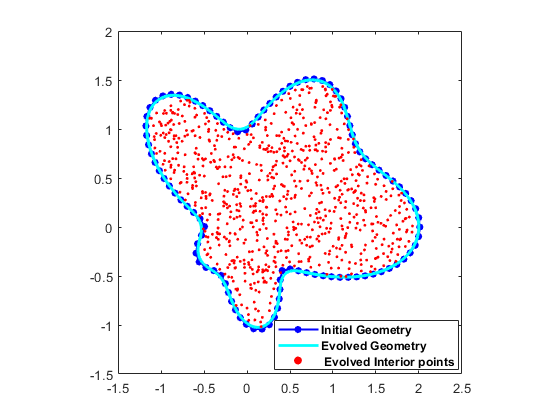}}\\
\subfigure[LOOCV and maximum principle error versus $r_s$]{
    \begin{overpic}[width=3in]{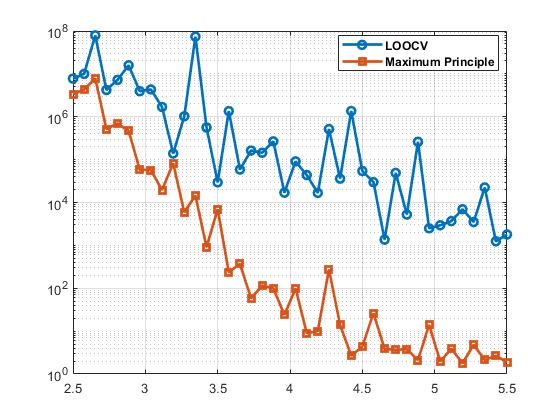}
        \put(20,240){\rotatebox{90}{\small $\ell_\infty$ norm error}}
        \put(500,10){$r_s$}
    \end{overpic}
}
\hspace{0.02in}
\subfigure[$E_{\mathrm{LOO},\infty}$ versus $r_s$ for various $h$]{
    \begin{overpic}[width=3in]{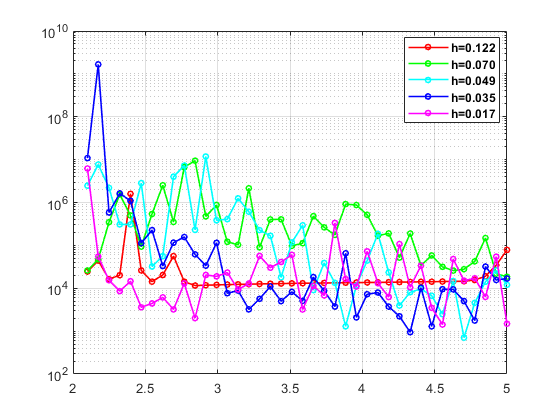}
        \put(20,220){\rotatebox{90}{\small $\ell_\infty$ norm error}}
        \put(500,10){$r_s$}
    \end{overpic}
}
\caption[Amoeba shape geometry with circular source points]{Amoeba shape geometry with circular source points
(a) Initial geometry with uniform distribution of interior points and source points (on a concentric circle)
(b) Evolved boundary and interior points at $t = T$
(c) LOOCV indicator $E_{\mathrm{LOO},\infty}$ and maximum principle error $E_{\mathrm{MP}}$ versus $r_s$
(d) $E_{\mathrm{LOO},\infty}$ versus $r_s$ for various fill distances $h$.
}
\label{fig:example4}
\end{figure}

The impact of this irregularity is immediately evident in the error indicators. As shown in Figure~\ref{fig:example4}(c), the LOOCV error $E_{\mathrm{LOO},\infty}$ exhibits large and erratic fluctuations, consistently overestimating the maximum principle error by several orders of magnitude across all tested values of $r_s$. This marked disparity signals a fundamental breakdown in the reliability of LOOCV as an error estimator for the square MFS system in such complex geometries. When this behavior emerges, it serves as a clear diagnostic that the standard formulation is no longer adequate, and that transitioning to an overdetermined system is necessary to restore numerical stability and trustworthy error quantification.

\subsubsection*{Summary}
Throughout this section, we have demonstrated that the LOOCV indicator remains a robust and practical tool for error estimation and parameter selection in the MFS, even as the domain geometry departs from classical settings. For smooth and mildly perturbed boundaries, LOOCV reliably tracks the true maximum principle error and provides consistent guidance across a range of parameters. As the geometry becomes increasingly non-smooth or non-convex, LOOCV continues to serve as an effective upper bound, so long as the interior discretization remains well distributed. However, in extreme cases where point clustering and irregularity arise, such as in the amoeba domain, the LOOCV can become overly pessimistic and lose its predictive value, signaling the necessity of switching to an overdetermined formulation. These results collectively establish LOOCV as a valuable and diagnostic indicator for the practical deployment of MFS in a broad spectrum of moving boundary problems.

\subsection{Improving MFS Robustness in Complex Domains via Zero-Padded Square Systems}
In the preceding square MFS experiments ($N_c = N_s$), exact collocation produced satisfactory interior-point distributions for the circle and the smooth/sharp-feature asterisks. However, for the amoeba-like domain, the same square formulation led to progressively scattered interior points. This can be explained using an effective-condition-number ($\kappa_{\mathrm{eff}}$) viewpoint. When the prescribed boundary data are the trace (and, for Neumann, the normal derivative) of a function that is harmonic in a neighborhood outside $\Omega$ (i.e. it admits a harmonic extension outside the domain), $\kappa_{\mathrm{eff}}$ is large: the data align mainly with well-resolved singular directions, so exact interpolation can reach near machine precision even when $\kappa(A)$ is large (see \cite{drombosky2009mfs,wong2011optimality}). By contrast, curvature-driven velocities and their numerical extensions typically introduce non-harmonic content into the right-hand side; this content projects onto left singular vectors associated with small singular values. The exact pointwise collocation then amplifies these components, injecting oscillations into $u$ and causing nonuniform transport of the interior points. In such cases, relaxing pointwise enforcement and solving in the least-squares sense, either by overcollocation ($N_c > N_s$) or, minimally, by augmenting the square system with a single zero row to induce least-squares behavior while keeping $N_c = N_s$, suppresses small-singular (oscillatory) modes, yields a smoother $u$, and maintains well-distributed interior points during evolution.

A prominent point in this section is the continued applicability of LOOCV–$\ell_\infty$ in an overdetermined setting (zero‑padded). In square systems, we report LOOCV–$\ell_\infty$ using the Rippa formula alongside the maximum‑principle indicator. However, once we move to overdetermined / least squares (including zero-row augmentation), the standard Rippa algorithm is not applicable. Instead, we can compute LOOCV using the hat‑matrix formulation, which produces an a posteriori error indicator that can be monitored throughout the evolution of complex geometries. As a low‑cost alternative, we also tested a modified pinv–Rippa algorithm designed to mimic Rippa in this setting; as shown in Figure~\ref{fig:example5over}(c), LOOCV–$\ell_\infty$ through the hat matrix provides a consistent bound throughout $r_s$, while pinv-Rippa is erratic and non-predictive.

\subsubsection*{Example 5: Amoeba Domain with Zero-Padded Square System}
\begin{figure}[t!]
\subfigure[Initial domain]{\includegraphics[width=3in]{fig/amoeba_squareCircle11.png}}
\subfigure[Evolved domain at $t=T$]{\includegraphics[width=3in]{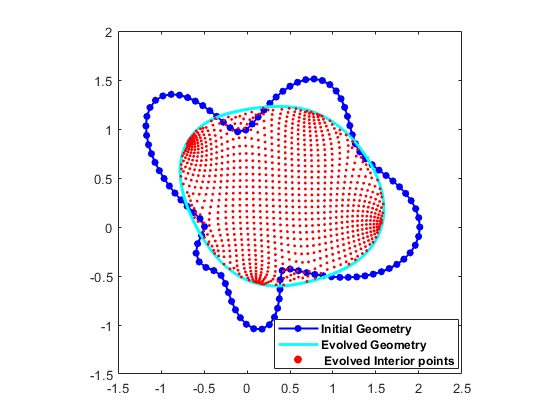}}
\begin{center}
\subfigure[LOOCV, pinv-Rippa and maximum principle error versus $r_s$]{
    \begin{overpic}[width=2.7in]{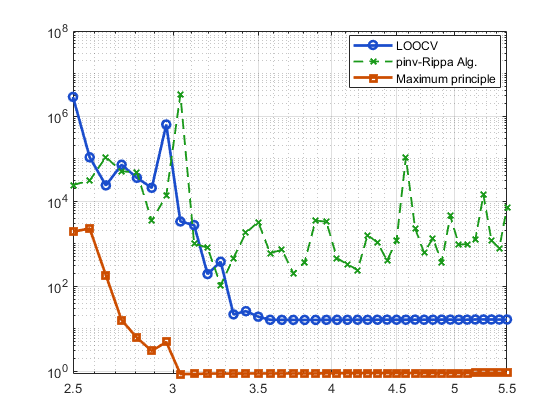}
        \put(20,240){\rotatebox{90}{\small $\ell_\infty$ norm error}}
        \put(500,10){$r_s$}
    \end{overpic}
}
\hspace{0.02in}
\subfigure[$E_{\mathrm{LOO},\infty}$ versus $r_s$ for various $h$]{
    \begin{overpic}[width=2.7in]{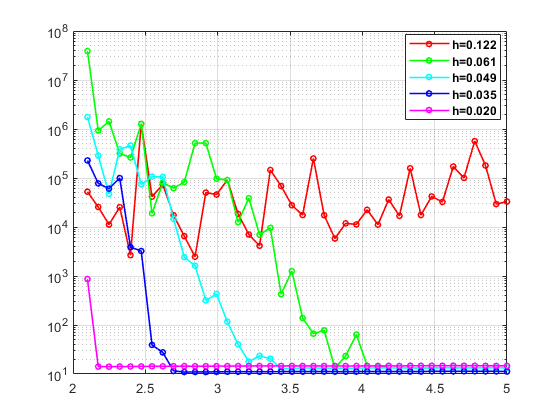}
        \put(20,220){\rotatebox{90}{\small $\ell_\infty$ norm error}}
        \put(500,10){$r_s$}
    \end{overpic}
}
\end{center}
\caption[Amoeba geometry with overdetermined system]{Amoeba geometry with overdetermined system
(a) Initial geometry with uniform distribution of interior points and circular source points of radius $r_s > 1$
(b) Evolved domain at $t = T$ with  improved well distributed interior points
(c) maximum principle error $E_{\mathrm{MP}}$ versus $r_s$.
}
\label{fig:example5over}
\end{figure}

Using the zero-padded MFS formulation for the amoeba-shaped domain leads to a marked improvement in the computed point distribution. The interior points, which previously showed scattering and clustering with the square formulation, now remain well distributed during evolution; see Figure~\ref{fig:example5over}(b). This follows from inducing a least-squares solve via zero padding while keeping the source and collocation sets unchanged.

Figure~\ref{fig:example5over}(c) plots the LOOCV indicator $\ell_\infty$, the pinv--Rippa variant, and the maximum-principle error versus the source radius $r_s$. In the zero-padded setting, the LOOCV indicator, computed from the hat matrix, provides an a posteriori measure that tracks changes in boundary error as the geometry evolves. For reference, we also report a pseudoinverse-based pinv--Rippa value; however, the classical Rippa identity does not apply to zero-padded least-squares systems, so this variant lacks theoretical justification and behaves inconsistently. Over a wide range of $r_s$, the LOOCV indicator remains conservative relative to the maximum-principle error, while the pinv--Rippa values fluctuate without reliable correspondence. These observations support the use of LOOCV as the primary a posteriori indicator for parameter selection and error control in complex, evolving domains.

\begin{remark}
Fixed circular source placement is often adequate for domains with simple or slowly changing boundaries. When the boundary becomes intricate and is represented by a dynamically adapted B-spline, requiring frequent insertion or deletion of boundary nodes to maintain mesh quality, the limitations of a fixed auxiliary circle become apparent. In such cases, fixed source locations may under-resolve regions of high curvature or rapid geometric change, leading to poor coverage of the boundary data. A practical remedy is to switch to geometry-adapted source placement, such as normal-offset families or related strategies (cf.\ \cite{chen2016choosing}). If these options remain insufficient for highly dynamic or complex domains, it can be preferable to replace the MFS with a fully ALE-based finite element method (e.g., ALE-FEM).
\end{remark}

\section{Numerical comparison with a classical ALE–FEM mesh motion}\label{sec:ale_fem_comparison}

The numerical experiments in Section~\ref{sec:nu_ex_5} showed that, once the zero–padded square formulation is adopted, the B-spline–based ALE–MFS produces a smooth harmonic-extension velocity and maintains a well-distributed set of interior points even for highly non-convex domains. In this section we use this mesh–motion strategy as a building block for a direct comparison with a classical, BGN–based ALE–FEM. The aim is to assess, on identical initial meshes, how the two approaches differ in terms of mesh quality during curvature–driven evolution.
\subsection{Mesh motion generated by BGN-based ALE--FEM}
\label{subsec:ale_fem_bgn}
To benchmark our ALE--MFS strategy, we compare against a classical ALE--FEM. The boundary evolution is treated by a BGN-type parametric FEM enforcing curvature-driven motion~\cite{barrett2007parametric,barrett2008parametric}, while the interior mesh velocity is computed as the harmonic extension of the boundary velocity and discretized with an evolving FEM. This results in a standard ALE mesh update that transports all mesh vertices accordingly. 

Let $\Omega(t)\subset\mathbb{R}^2$ be a smoothly evolving domain with boundary $\Gamma(t)=\partial\Omega(t)$ for $t\in[0,T]$. The interior (mesh) velocity $\vec u(\cdot,t):\Omega(t)\to\mathbb{R}^2$ is defined as the harmonic extension of the boundary velocity $\vec V_\Gamma(\cdot,t)$, as given in equation \eqref{eq:harmonic-extension}. Let $\mathcal{A}_t:\Omega^0\times[0,T]\to\mathbb{R}^2$ be the ALE flow map governed by velocity $\vec u$. Then $\mathcal{A}_t$ maps $\Omega^0$ to $\Omega(t)$ and $\Gamma^0$ to $\Gamma(t)$, respectively, and satisfies the \eqref{eq:ALEframework}. Then, the ALE weak formulation for such harmonic extension is defined as

\begin{equation}
   \begin{cases} 
\frac{\partial}{\partial t} A_t = \vec{u} \circ A_t  \\ 
(\nabla \vec{u}, \nabla \vec{x}) = 0 & \forall x \in H_0^1(\Omega(t)) \\ 
\vec{u} = \vec{v}_\Gamma & \text{on } \Gamma(t)
\end{cases} 
\label{eq:weakFormHarmonicExtension}
\end{equation}
with initial conditions \( A_t(0) = \vec{id}|_{\Omega^0} \) and $\mathrm{id}(x)=x$ denotes the identity on $\mathbb{R}^2$. We will use the evolving ALE-FEM \cite{li2025optimal,gong2024convergent} to compute the discrete velocity $u_h^n$ based on the ALE weak formulation \eqref{eq:weakFormHarmonicExtension}.

The boundary motion itself is obtained from the curvature-driven BGN-type parametric FEM proposed by Barrett, Garcke, and N\"urnberg ~\cite{barrett2007parametric,barrett2008parametric}. In this approach, the evolving interface $\Gamma(t)$ is approximated by a polygonal curve $\Gamma_h(t)$, whose vertices are updated at each time step by solving a discrete variational problem that weakly enforces the curvature law and introduces a tangential redistribution of nodes. At the continuous level, the curvature flow is described by
\begin{equation}
    \begin{cases}
        \vec{V}_\Gamma\cdot\vec{n} = -\,\kappa & \text{on } \Gamma(t),\\
        \kappa\,\vec{n} = -\,\Delta_{\Gamma(t)} \,\mathrm{id} & \text{on } \Gamma(t),
    \end{cases}
\label{eq:BGNFormulation}
\end{equation}
where $\kappa$ and $\vec{n}$ are the scalar curvature and outward unit normal on $\Gamma(t)$, respectively, and $\Delta_{\Gamma(t)}$ is the Laplace--Beltrami operator along $\Gamma(t)$. The discrete BGN scheme is known to generate a nearly equidistributed set of nodes along the curve and is widely used as a robust interface approximation for geometric flows, see~\cite{bao2021structure}.

To transfer this boundary evolution to a volume mesh, we consider an initial polygonal domain $\Omega_0\subset\mathbb{R}^2$ with boundary $\Gamma_0$ and a conforming, shape--regular triangulation $\mathcal{T}_h^0$ of $\Omega_0$. The triangles of $\mathcal{T}_h^0$ satisfy a minimal angle condition $\theta_{\min}(\mathcal{T}_h^0) \geq 20^\circ,$ where $\theta_{\min}(\mathcal{T}_h)$ denotes the smallest interior angle over all triangles in $\mathcal{T}_h$. At time $t_n$, the BGN scheme supplies an updated polygonal curve $\Gamma_h^n$ together with nodal velocities $\vec V_{\Gamma,h}^n$ at the curve vertices. The interior mesh velocity $\vec u_h^n$ is then defined as the solution of the discrete harmonic extension problem
\begin{equation}
\begin{cases}
(\nabla \vec{u}_h^n, \nabla \xi_h) = 0 &\ \forall \xi_h \in \mathcal{S}_h [\vec{x}_h^n]\\[0.25em]
\vec u_h^n = \vec V_{\Gamma,h}^n &\text{on } \Gamma_h^n,
\end{cases}
\label{eq:ale_fem_harmonic}
\end{equation}
where $\mathcal{S}_h [\vec{x}_h^n]$ is the finite element space define on the evolving domain $\Omega_h^n$. The mesh vertices are then updated by the time--stepping rule
\begin{equation}
  \vec x_h^{n+1} = \vec x_h^n + \Delta t\, \vec u_h^n(\vec x_h^n),
  \label{eq:ale_fem_update}
\end{equation}
which defines the ALE--FEM mesh--motion operator used as benchmark in the following.

\subsection{B-spline based ALE--MFS mesh motion on a given FEM mesh}

For a consistent comparison, both the BGN-based ALE--FEM method and the proposed B-spline based ALE--MFS method are started from the same initial triangulations $\mathcal{T}_h^0$ of the reference domain. We consider two representative mesh sizes, $h=0.2$ and $h=0.05$. In each case the boundary is given by the vertices on $\Gamma_0$, and both methods are driven by the same curvature--flow law \eqref{eq:BGNFormulation} in the normal direction.

In the B-spline based  ALE--MFS approach, the boundary at time $t_n$ is described as a point cloud. For each boundary point $\vec q_i\in\Gamma_h^n$ the curvature $\kappa(\vec q_i,t_n)$ and outward unit normal $\vec n(\vec q_i,t_n)$ are obtained from the local B-spline reconstruction detailed in Section~\ref{subsec:mfs_unprescribed} and \cite{ammad2025eabe}. The boundary velocity is then prescribed as
\begin{equation}
  \vec V_\Gamma^{\mathrm{MFS}}(\vec q_i,t_n)
  =
  -\,\kappa(\vec q_i,t_n)\,\vec n(\vec q_i,t_n),
  \quad \vec q_i\in\Gamma_h^n.
\end{equation}
The associated interior velocity field $\vec u^{\mathrm{MFS}}$ is defined as the solution of
\begin{equation}
\begin{cases}
\Delta \vec u^{\mathrm{MFS}} = 0 &\text{in } \Omega_h^n,\\[0.25em]
\vec u^{\mathrm{MFS}} = \vec V_\Gamma^{\mathrm{MFS}} &\text{on } \Gamma_h^n,
\end{cases}
\end{equation}
and is approximated by the MFS ansatz \eqref{eq:mfs-ansatz} with sources placed on a fixed exterior circle as in Section~\ref{sec:points_mfs}. The coefficients of the MFS expansion are computed using the zero--padded square formulation introduced in Section~\ref{sec:error}, which yields a smooth, least--squares harmonic extension. Once $\vec u^{\mathrm{MFS}}$ is available, the nodal positions of the mesh are updated using the same ALE update rule \eqref{eq:ale_fem_update}, with $\vec u_h^n$ replaced by $\vec u^{\mathrm{MFS}}(\cdot,t_n)$ evaluated at the mesh vertices.

At $t=0$, both methods share the same triangulation and boundary. From that point on, one method uses BGN-based ALE-FEM to propagate the mesh, while the other uses the B-spline based ALE--MFS mesh motion. This setting isolates the effect of the mesh--motion strategy on mesh quality.

\subsection{Mesh--quality indicators}

To compare the quality of the evolving meshes, we monitor two standard indicators throughout the numerical evolution. Let $\mathcal{T}_h^n$ be the triangulation at time $t_n$.

First, we consider the minimum interior angle
\[
\theta_{\min}^n = \min_{K\in\mathcal{T}_h^n}\ \min_{1\leq\ell\leq 3} \theta_\ell(K),
\]
where $\theta_\ell(K)$ is the $\ell$-th interior angle of the triangle $K$. Larger values of $\theta_{\min}^n$ correspond to better-shaped triangles, whereas very small angles signal deteriorated elements. A lower bound on $\theta_{\min}^n$ is a classical requirement for stability and optimal convergence of finite element discretizations; see, e.g., the angle and shape-regularity conditions in \cite{li2025optimal,duan2022energy}. In the figures we therefore display the reference line $\theta_{\mathrm{thr}} = \pi/18 \approx 10^\circ$, and interpret the crossing of this threshold as an indication that remeshing would be required in a standard ALE--FEM pipeline.

Second, we use the mesh ratio
\[
\rho_{\mathrm{mesh}}^n
=
\frac{h_{\max}^n}{h_{\min}^n},
\qquad
h_{\max}^n = \max_{K\in\mathcal{T}_h^n}\mathrm{diam}(K),\quad
h_{\min}^n = \min_{K\in\mathcal{T}_h^n}\mathrm{diam}(K),
\]
which measures the spread of element sizes in the mesh. Values close to one indicate quasi--uniform meshes and are therefore preferable, while large values point to the coexistence of extremely small and comparatively large elements, a situation that typically harms conditioning and accuracy.

Both indicators are reported as functions of time for the two geometries (sharp asterisk and amoeba) and the two mesh sizes $h=0.2$ and $h=0.05$.

\subsection{Results for the sharp asterisk}

\begin{figure}[H]
  \centering
  \subfigure[Initial geometry and mesh at $t=0$ for $h=0.2$.]{
    \includegraphics[width=3in]{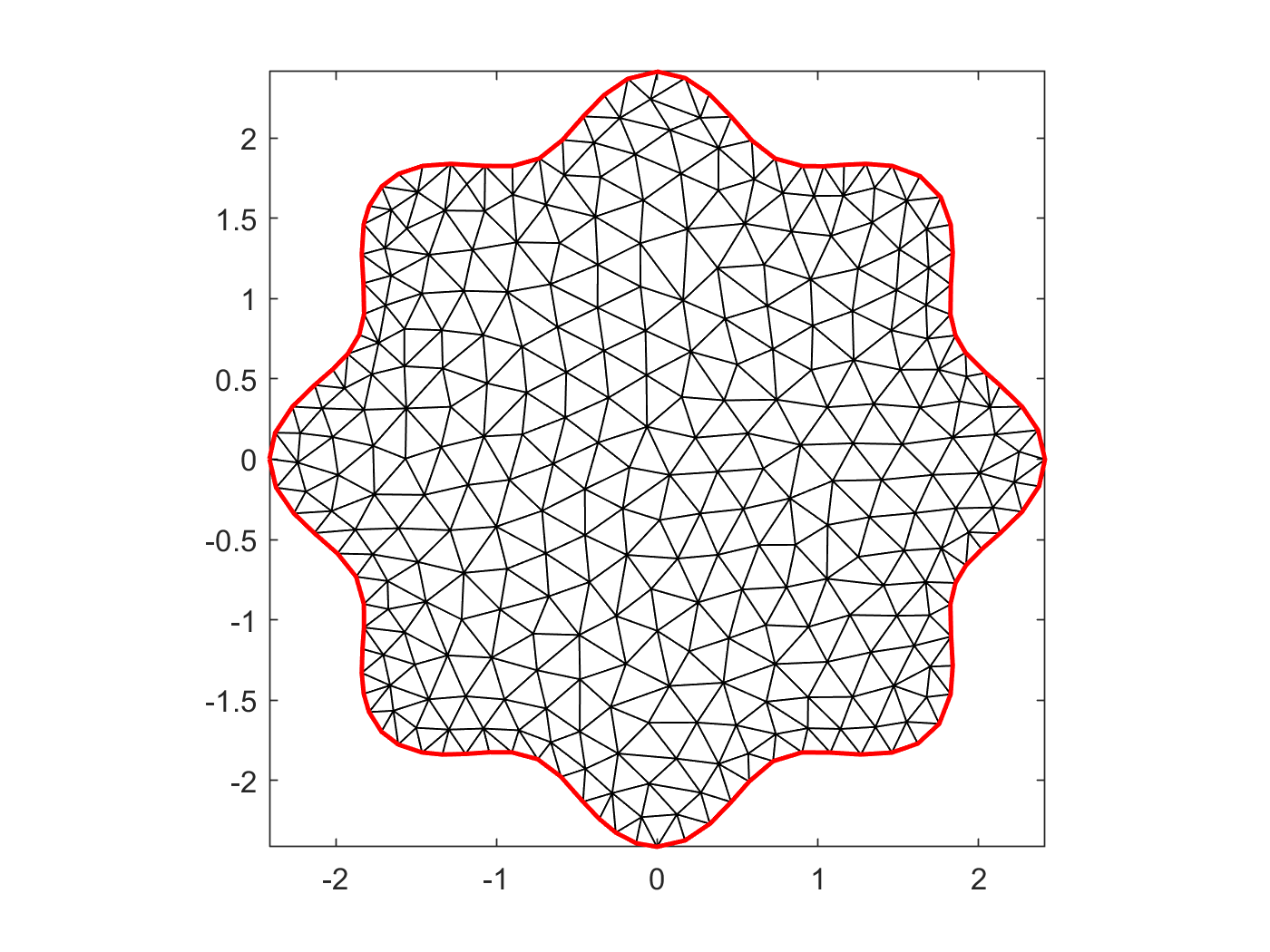}
  }
  \subfigure[Evolved geometry and mesh at $t=T$ for $h=0.2$.]{
    \includegraphics[width=3in]{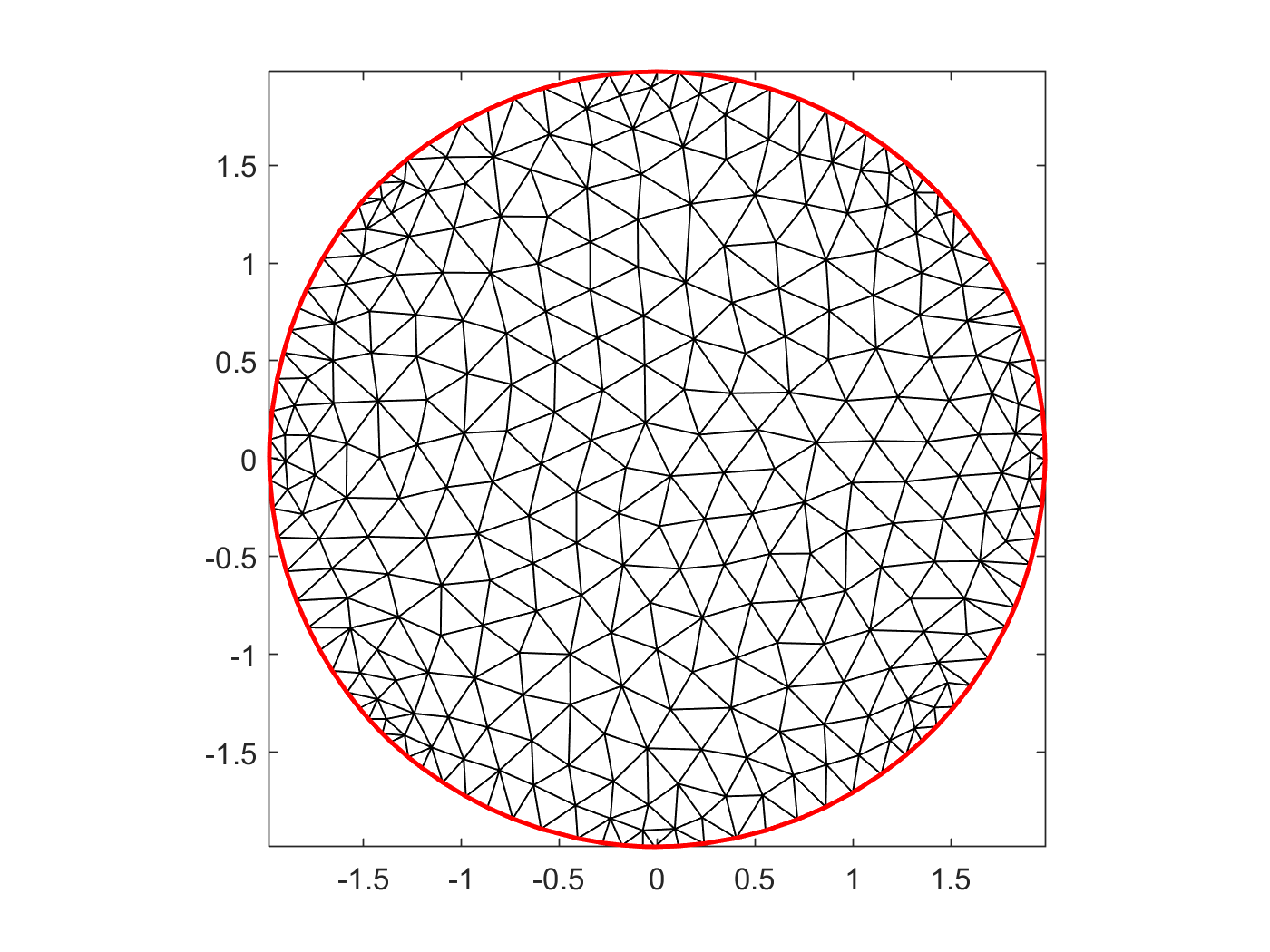}
  }
  \subfigure[Minimum triangle angle for $h=0.2$.]{
    \begin{overpic}[width=3in]{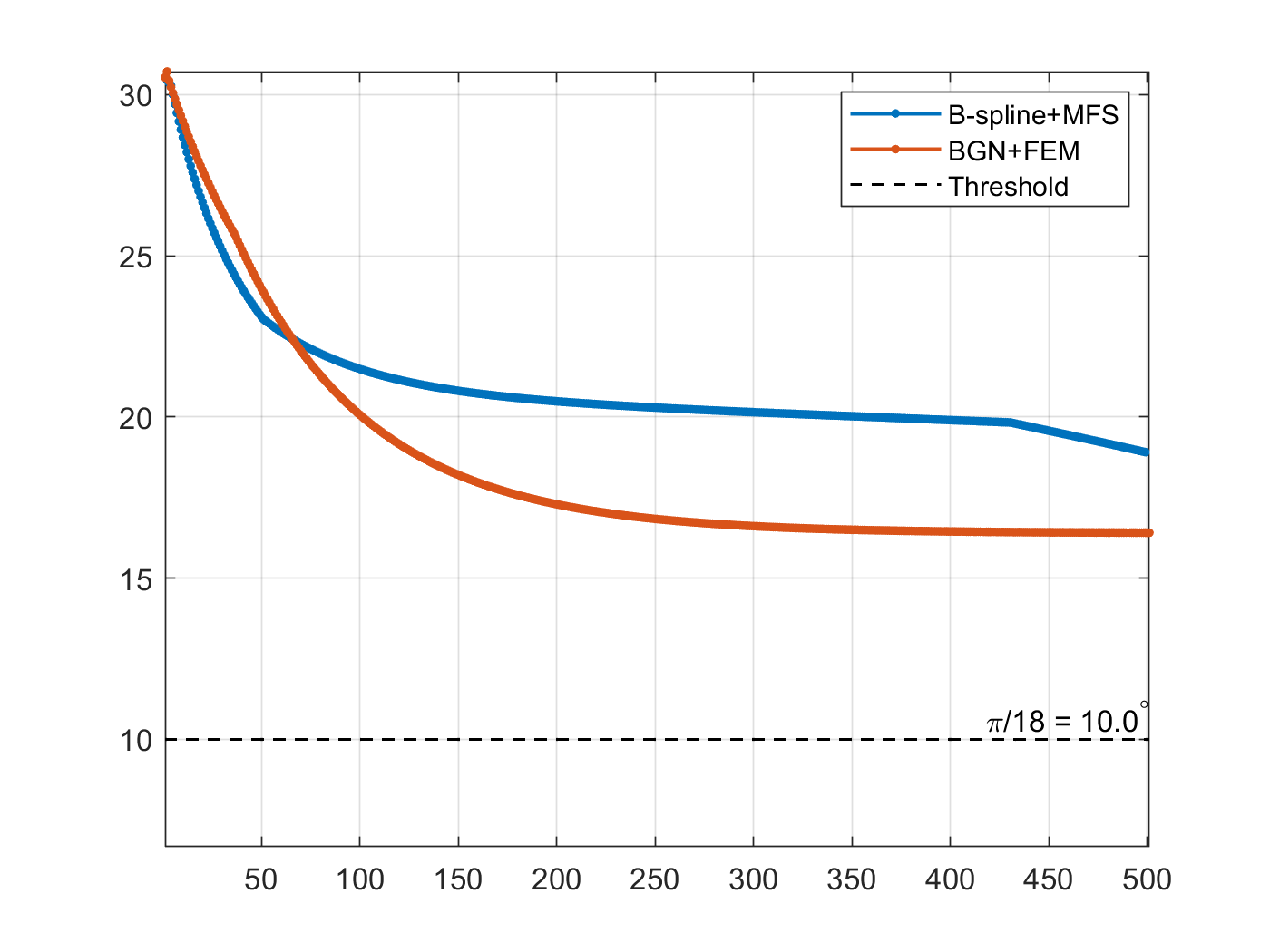}
      % x-axis label
      \put(480, -3){\scriptsize Iteration}
      % y-axis label (rotated)
      \put(50, 300){\rotatebox{90}{\scriptsize Angle (deg)}}
    \end{overpic}
  }
  \subfigure[Minimum triangle angle for $h=0.05$.]{
    \begin{overpic}[width=3in]{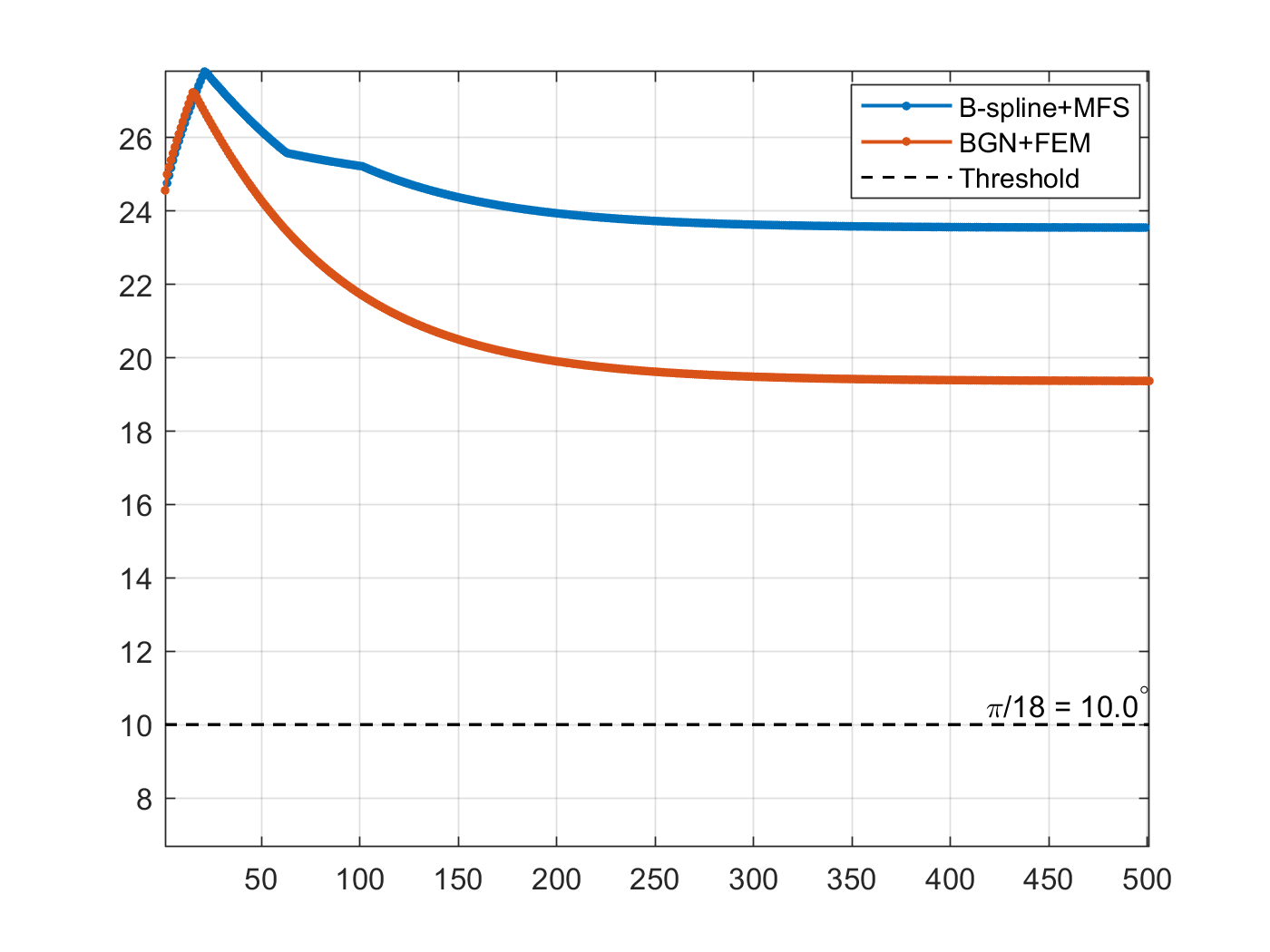}
      % x-axis label
      \put(480, -3){\scriptsize Iteration}
      % y-axis label (rotated)
      \put(50, 300){\rotatebox{90}{\scriptsize Angle (deg)}}
    \end{overpic}
  }
  \subfigure[Mesh ratio for $h=0.2$.]{
    \begin{overpic}[width=3in]{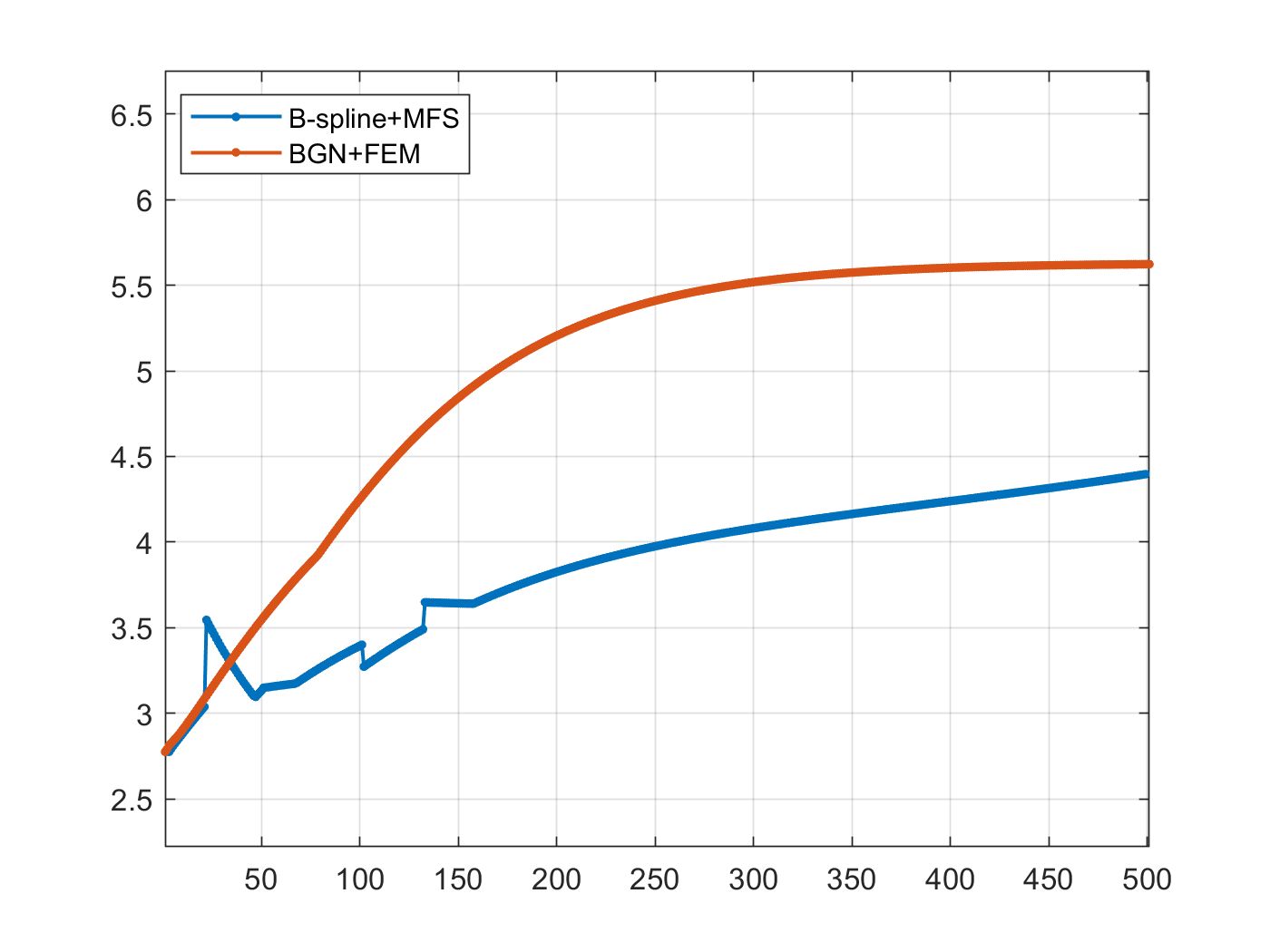}
      % x-axis label
      \put(480, -3){\scriptsize Iteration}
      % y-axis label (rotated)
      \put(50, 300){\rotatebox{90}{\scriptsize Mesh ratio}}
    \end{overpic}
  }
  \subfigure[Mesh ratio for $h=0.05$.]{
    \begin{overpic}[width=3in]{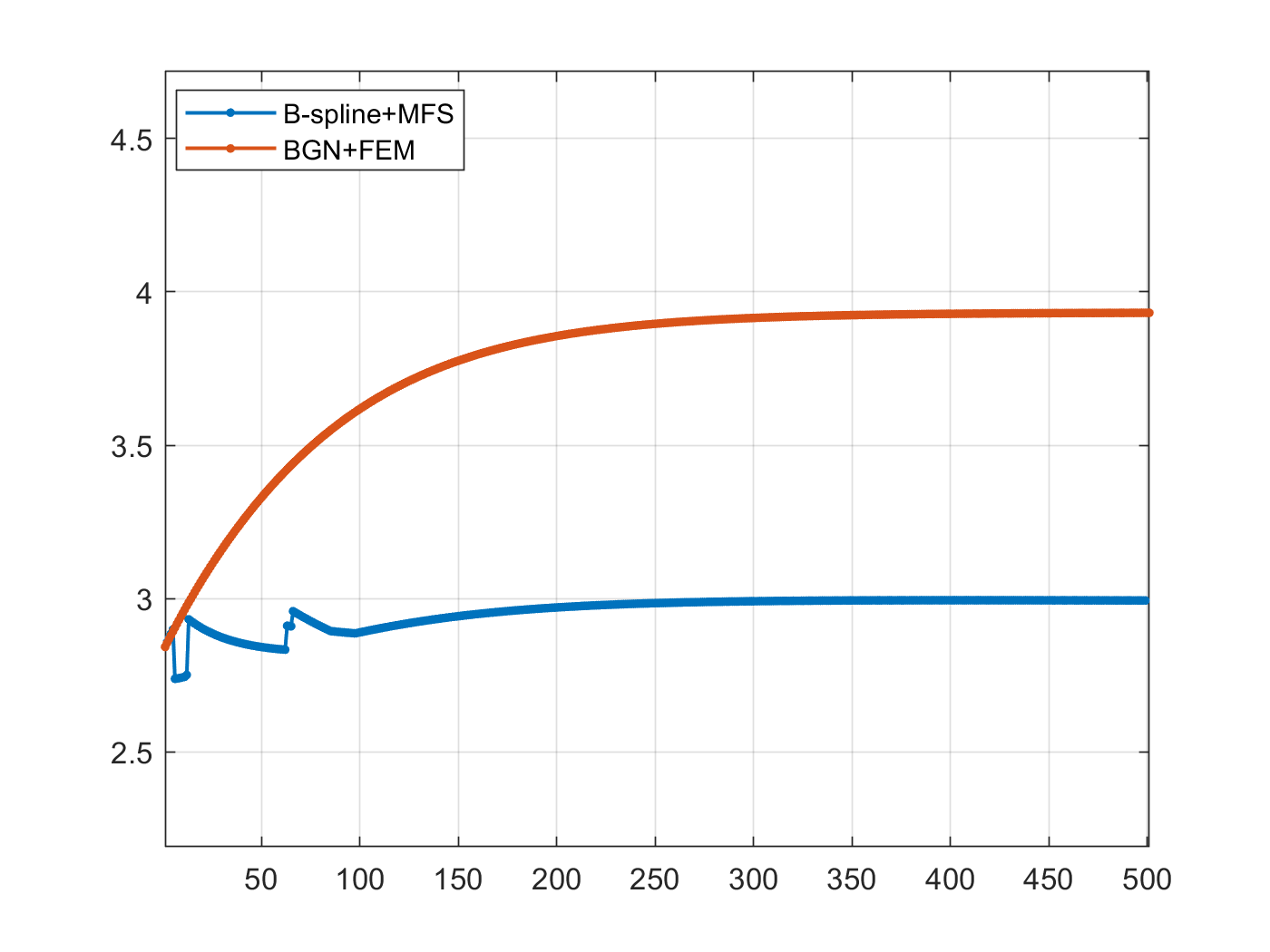}
      % x-axis label
      \put(480, -3){\scriptsize Iteration}
      % y-axis label (rotated)
      \put(50,300){\rotatebox{90}{\scriptsize Mesh ratio}}
    \end{overpic}
  }
\caption{Sharp asterisk geometry. Panels~(a)--(b): initial and evolved meshes for the coarse mesh size $h=0.2$. Panels~(c)--(d): evolution of the minimum triangle angle $\theta_{\min}^n$ (in degrees) versus iteration for $h=0.2$ and $h=0.05$. Panels~(e)--(f): evolution of the mesh ratio $\rho_{\mathrm{mesh}}^n$ versus iteration for $h=0.2$ and $h=0.05$. In all panels, larger values of the minimum angle and smaller values of the mesh ratio indicate better mesh quality.}
  \label{fig:asterisk_comparison}
\end{figure}

We first consider the sharp asterisk domain introduced in Example~3. 
Panels~(a)--(b) of Figures~\ref{fig:asterisk_comparison}(c)--(d) display, for the coarse mesh size $h=0.2$, 
the initial triangulation $\mathcal{T}_h^0$ together with the boundary curve, 
and the mesh at the final time $t=T$. For the finer mesh $h=0.05$ the 
individual elements are too small to be informative in a global view, 
so we do not show the corresponding triangulations and only report the 
mesh--quality indicators. The tips of the asterisk yield regions of high 
curvature where local mesh deterioration can easily occur.

Figures~\ref{fig:asterisk_comparison}(c)--(d) show the evolution of the minimum interior angle $\theta_{\min}^n$ for the two mesh sizes. For $h=0.2$, both methods start from a minimum angle near $30^\circ$. As time progresses and the geometry smooths under curvature flow, the minimum angle decreases, but remains uniformly above the $10^\circ$ threshold. The BGN-based ALE-FEM curve decreases somewhat more rapidly and levels off at a lower value than the B-spline based ALE--MFS curve. For $h=0.05$, the difference becomes more visible: the B-spline based ALE--MFS mesh maintains a clearly larger minimum angle over most of the simulation interval, and in particular stays farther from the remeshing threshold. This indicates that the MFS--based harmonic extension, combined with the B-spline boundary reconstruction, tends to avoid the formation of very acute triangles near the tips.

The corresponding mesh ratios are shown in Figures~\ref{fig:asterisk_comparison}(e)--(f). For the coarse mesh, the BGN-based ALE-FEM ratio increases steadily and approaches a value around $4.5$ at final time, whereas the B-spline based ALE--MFS ratio stabilizes closer to $3.5$. For the fine mesh, both methods yield larger mesh ratios, as expected for a finer resolution of sharp features; nevertheless, the B-spline ALE--MFS curve settles at a smaller asymptotic value than BGN-based ALE-FEM. So, for the sharp asterisk geometry the two methods have comparable robustness, with a slight advantage for the MFS--based scheme in terms of minimum angle and uniformity of element sizes.

\subsection{Results for the amoeba domain}

The second comparison concerns the amoeba geometry, which is strongly non-convex and contains several lobes of different scales. This domain is more demanding for mesh motion than the asterisk: as curvature flow proceeds, thin protrusions shrink, merge and disappear, and the interior velocity field must accommodate large shape changes.

For the coarse mesh size $h=0.2$, the initial and final meshes are displayed in Figures~10(a)--(b). As in the asterisk case, we omit the $h=0.05$ triangulations from the figures, since the elements are too small to allow a meaningful visual comparison on the full domain, and instead focus on the quantitative indicators. Qualitatively, the B-spline based ALE--MFS mesh at $t=T$ appears to contain fewer extremely thin elements, particularly in regions where lobes have collapsed. This impression is supported by the quantitative indicators.

\begin{figure}[H]
  \centering
  \subfigure[Initial geometry and mesh at $t=0$ for $h=0.2$.]{
    \includegraphics[width=3in]{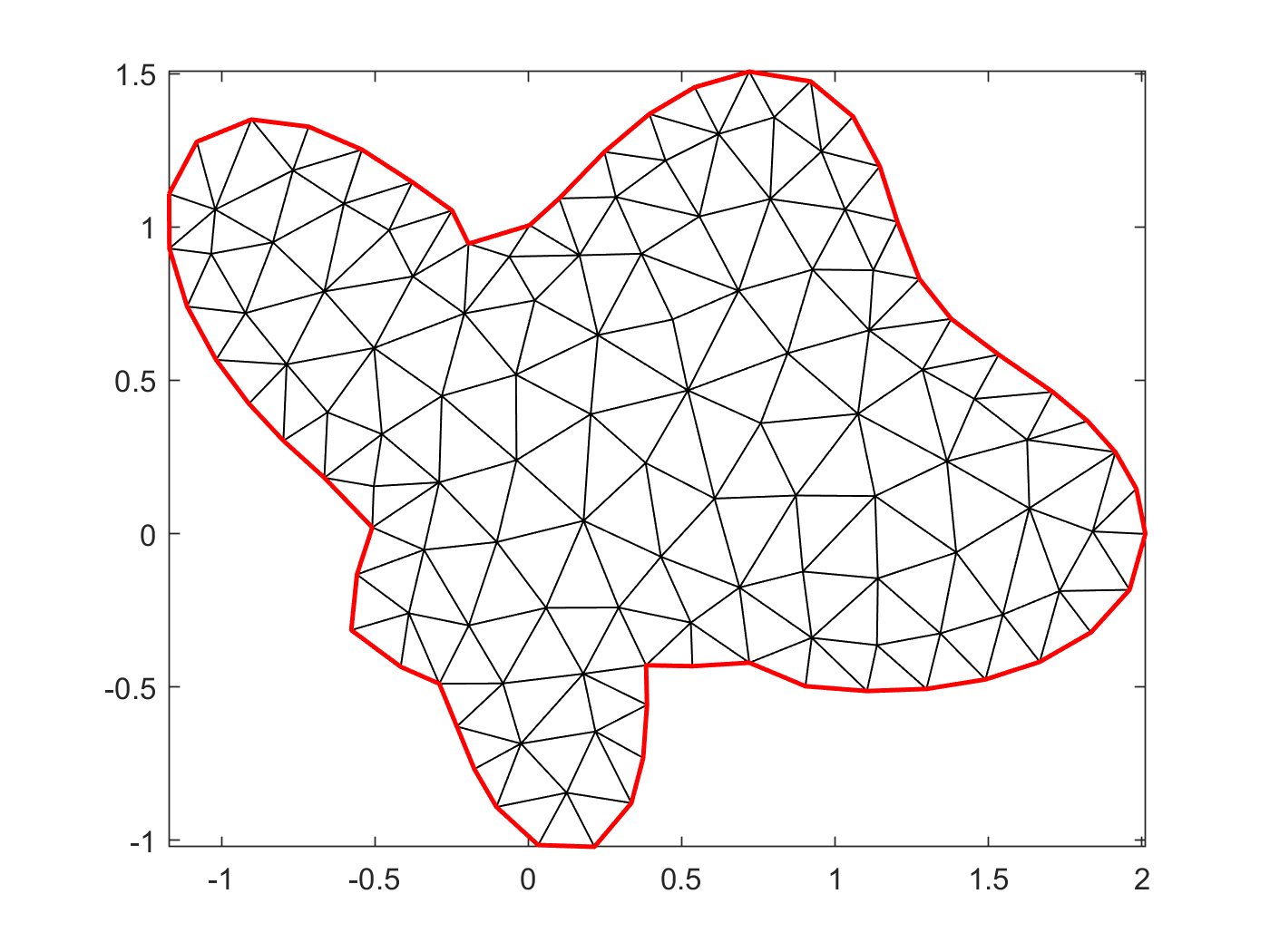}
  }
  \subfigure[Evolved geometry and mesh at $t=T$ for $h=0.2$.]{
    \includegraphics[width=3in]{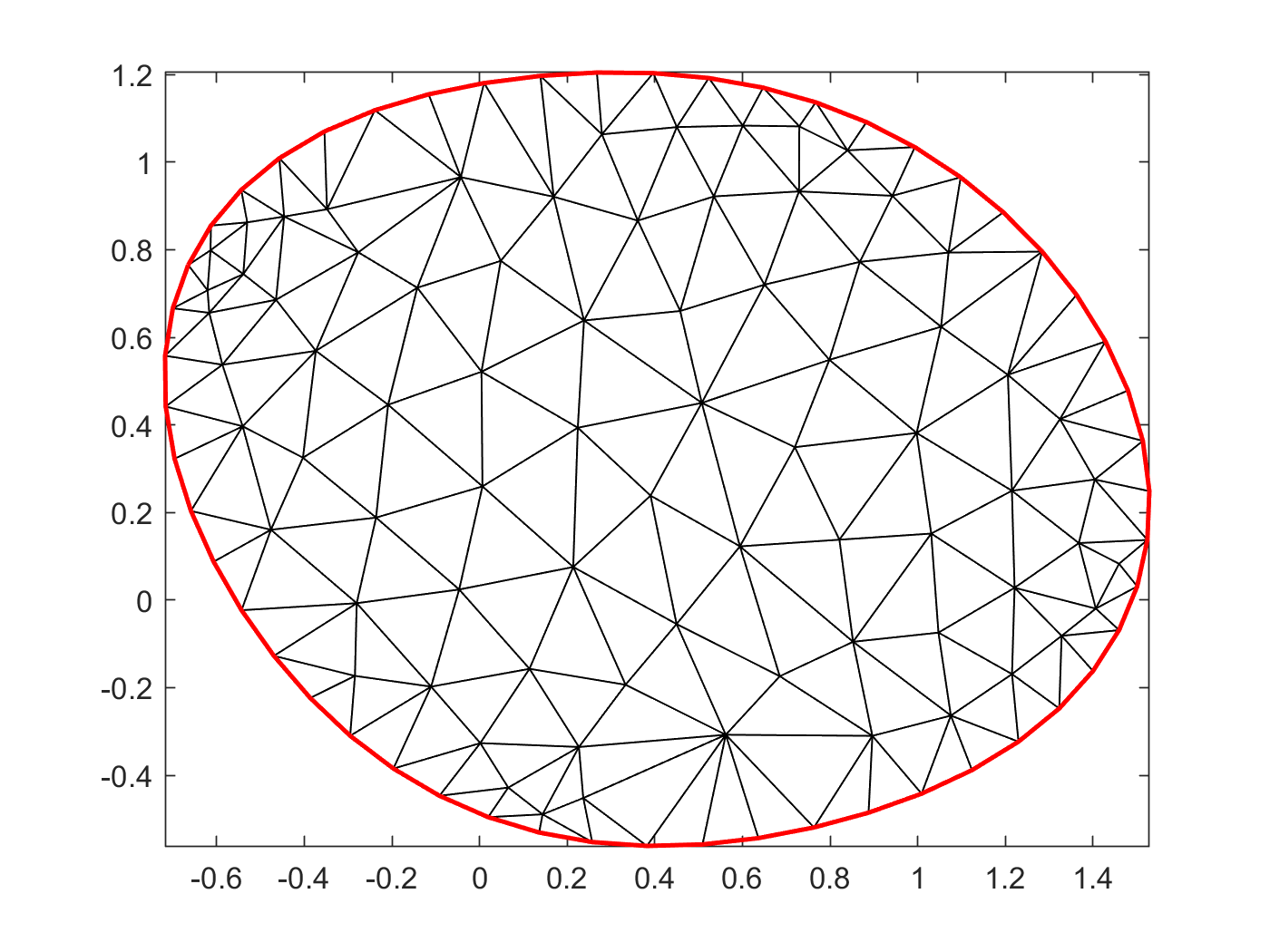}
  }
  \subfigure[Minimum triangle angle for $h=0.2$.]{
    \begin{overpic}[width=3in]{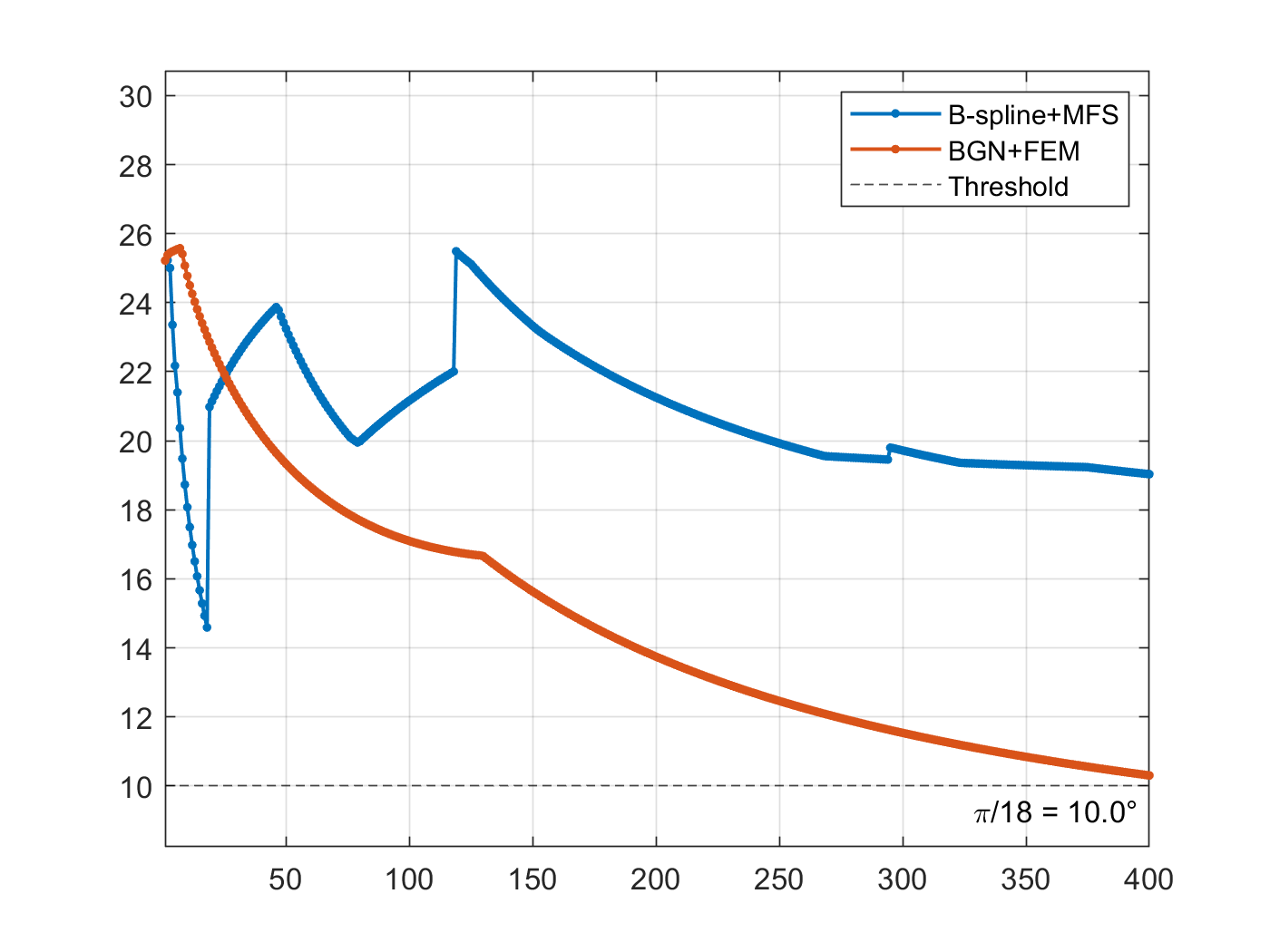}
      % x-axis label
      \put(480, -3){\scriptsize Iteration}
      % y-axis label (rotated)
      \put(50,300){\rotatebox{90}{\scriptsize Angle (deg)}}
    \end{overpic}
  }
  \subfigure[Minimum triangle angle for $h=0.05$.]{
    \begin{overpic}[width=3in]{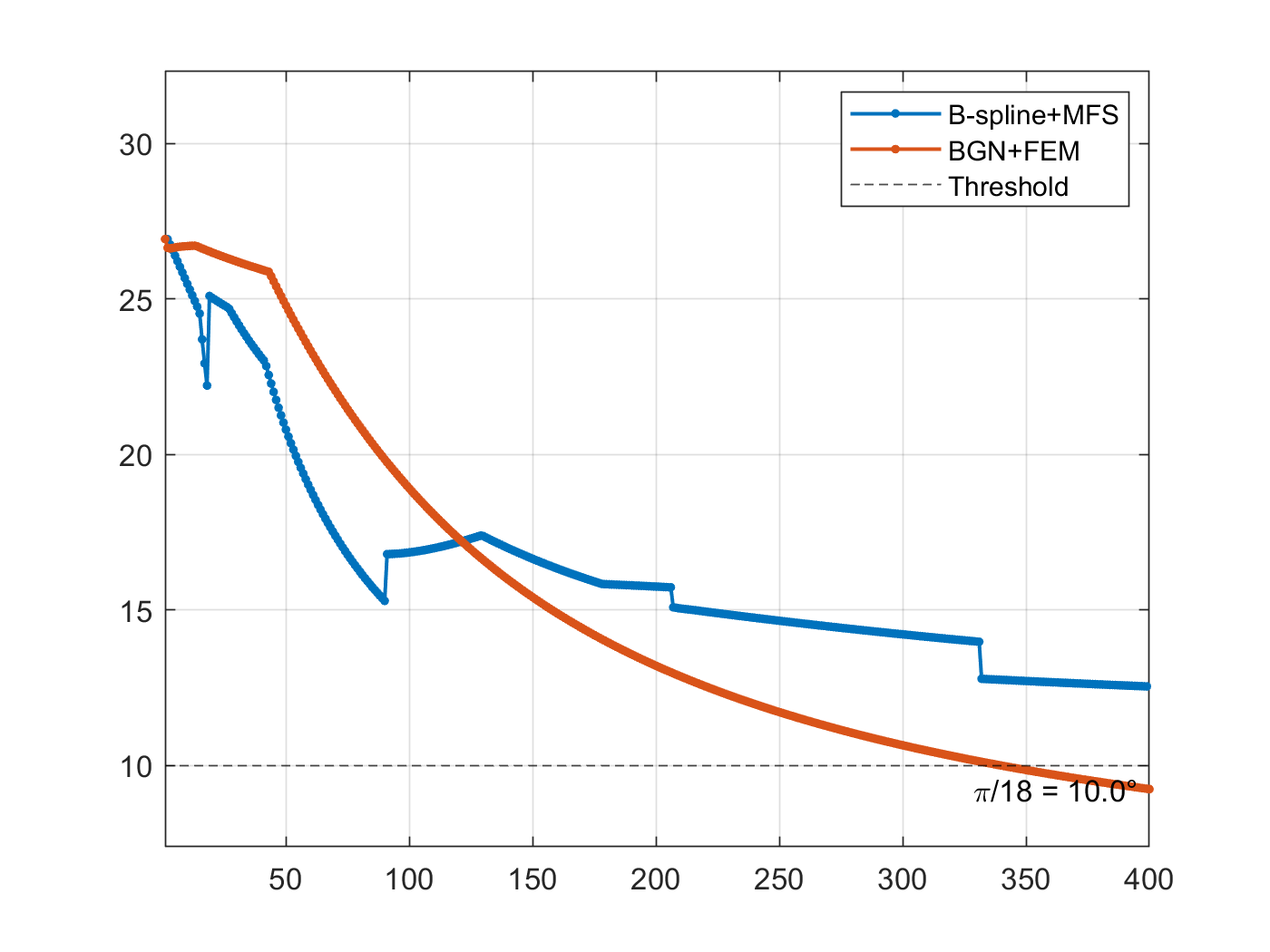}
      % x-axis label
      \put(480, -3){\scriptsize Iteration}
      % y-axis label (rotated)
      \put(50,300){\rotatebox{90}{\scriptsize Angle (deg)}}
    \end{overpic}
  }
  \subfigure[Mesh ratio for $h=0.2$.]{
    \begin{overpic}[width=3in]{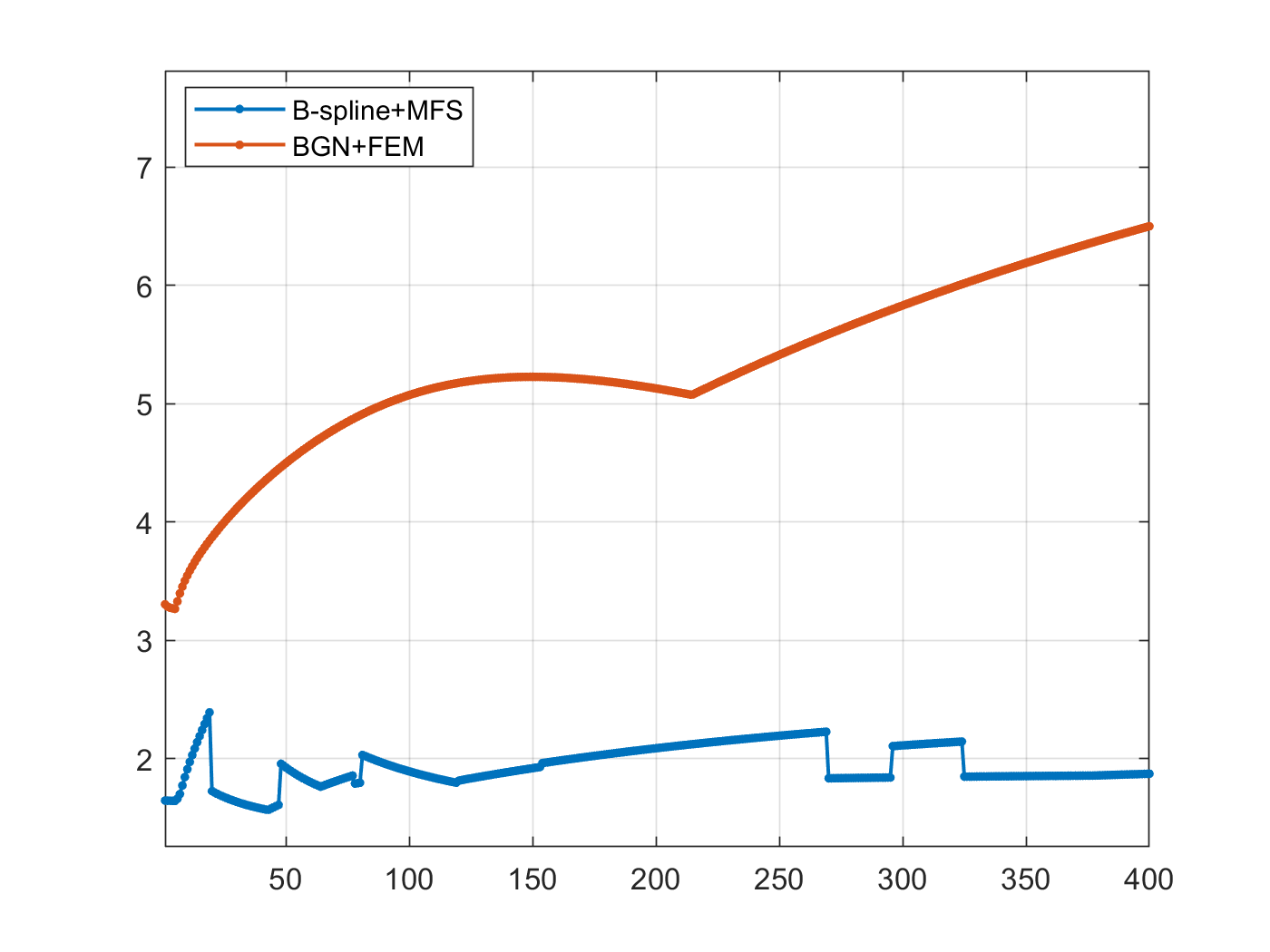}
      % x-axis label
      \put(480, -3){\scriptsize Iteration}
      % y-axis label (rotated)
      \put(50,300){\rotatebox{90}{\scriptsize Mesh ratio}}
    \end{overpic}
  }
  \subfigure[Mesh ratio for $h=0.05$.]{
    \begin{overpic}[width=3in]{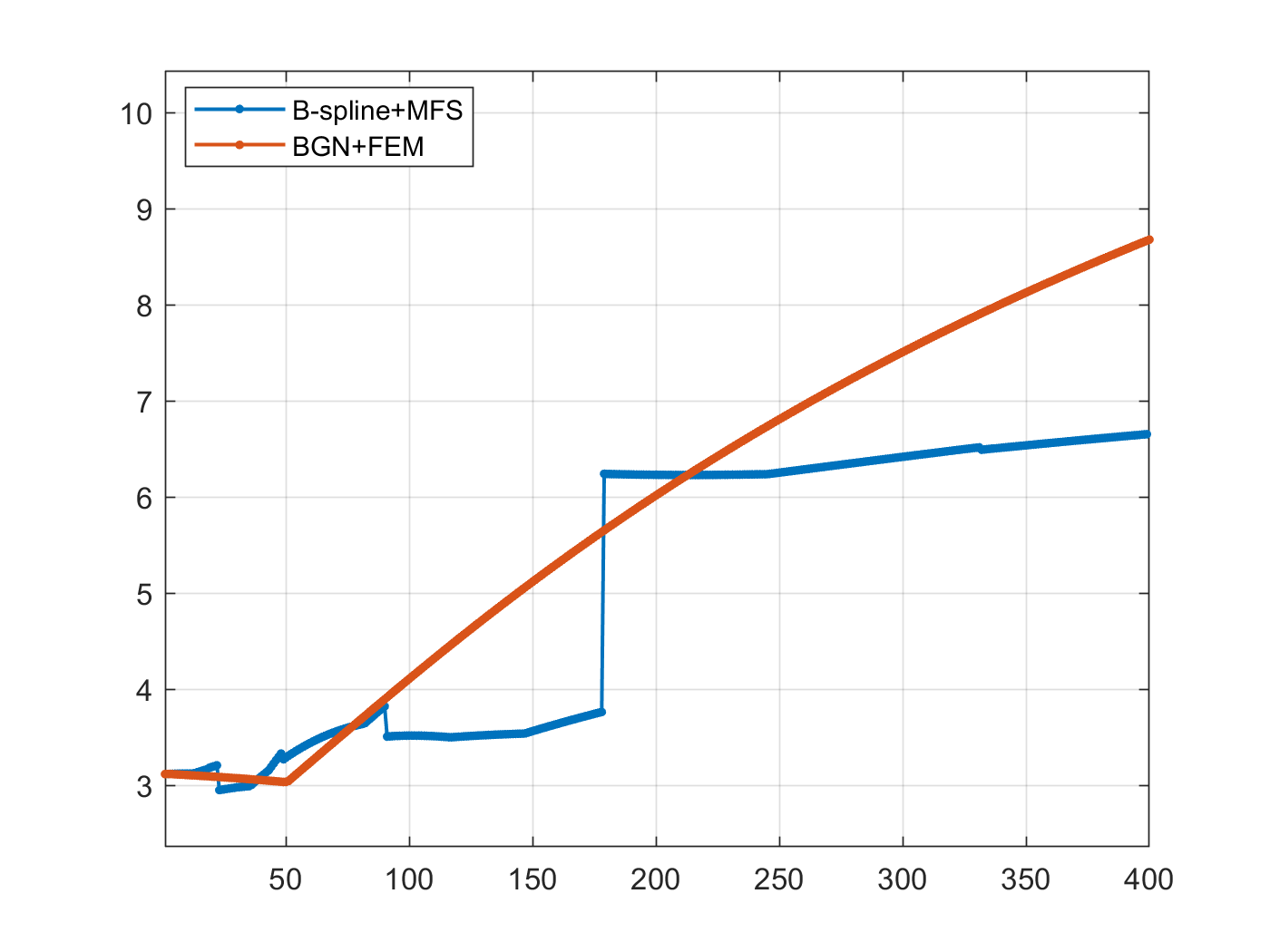}
      % x-axis label
      \put(480, -3){\scriptsize Iteration}
      % y-axis label (rotated)
      \put(50,300){\rotatebox{90}{\scriptsize Mesh ratio}}
    \end{overpic}
  }
\caption{Amoeba geometry. Panels~(a)--(b): initial and evolved meshes for the coarse mesh size $h=0.2$. Panels~(c)--(d): evolution of the minimum triangle angle $\theta_{\min}^n$ (in degrees) versus iteration for $h=0.2$ and $h=0.05$. Panels~(e)--(f): evolution of the mesh ratio $\rho_{\mathrm{mesh}}^n$ versus iteration for $h=0.2$ and $h=0.05$.}
\label{fig:amoeba_comparison}
\end{figure}

For $h=0.2$, the evolution of $\theta_{\min}^n$ is shown in Figures~\ref{fig:amoeba_comparison}(c). The BGN-based ALE-FEM minimum angle decreases almost monotonically and approaches the $10^\circ$ threshold as the domain becomes smoother. In contrast, the B-spline based ALE--MFS curve exhibits some oscillations but remains significantly above the threshold for the entire time interval. These oscillations are not a sign of instability; they are induced by the adaptive nature of the B-spline boundary reconstruction, which is allowed to insert or remove boundary points when they become too close or too far apart according to a fixed distance threshold. Whenever points are removed (or inserted), the local mesh connectivity is slightly rearranged, which produces small jumps in the minimum angle while still maintaining a globally good element quality. The effect is even more pronounced for the finer mesh $h=0.05$ in Figures~\ref{fig:amoeba_comparison}(d): while the BGN-based ALE-FEM minimum angle moves close to the $10^\circ$ line, the B-spline based ALE--MFS minimum angle stays appreciably larger over most of the evolution, indicating that the B-spline based ALE--MFS motion suppresses the formation of extremely acute elements even when long and thin features are present.

The mesh ratios in Figures~\ref{fig:amoeba_comparison}(e)--(f) highlight the same tendency. For $h=0.2$, the BGN-based ALE-FEM mesh ratio grows steadily and reaches values above $7$, signalling the appearance of very small elements adjacent to much larger ones. The B-spline based ALE--MFS mesh ratio, on the other hand, remains between $2.5$ and $3.5$ with only moderate variation. For $h=0.05$, the difference becomes substantial: the BGN-based ALE-FEM ratio continues to increase and exceeds $8$ towards the final time, whereas the B-spline based ALE--MFS ratio grows more slowly and stays markedly smaller. This behaviour suggests that the meshless harmonic extension generated by the MFS ansatz distributes the interior motion more evenly, so that elements do not become excessively stretched or compressed in localized regions.

\subsection{Discussion}
The comparison on these two geometries shows that, once a good interior point distribution is provided by the zero--padded formulation, the B-spline based  ALE--MFS mesh--motion strategy behaves at least as robustly as the classical BGN-based ALE-FEM approach, and in several cases yields visibly better meshes. In the moderately singular asterisk geometry both methods preserve acceptable element quality without remeshing, but the MFS--based velocity tends to keep larger minimum angles and smaller mesh ratios. In the more intricate amoeba geometry, the difference becomes more pronounced: the B-spline based ALE--MFS scheme maintains the minimum angle safely above the remeshing threshold and limits the growth of the mesh ratio, whereas the BGN-based ALE-FEM mesh enters a regime where many finite element codes would trigger remeshing based on standard quality criteria \cite{li2025optimal,duan2022energy}.

A further point is that the B-spline based ALE--MFS mesh motion can be integrated into an existing finite element code without altering the PDE discretization. The only change is the replacement of the FEM harmonic extension \eqref{eq:ale_fem_harmonic} by the meshless harmonic extension obtained from the MFS ansatz. In this sense, the proposed B-spline based ALE--MFS method functions as a drop--in alternative mesh--motion module, with the potential to postpone or even avoid expensive remeshing steps in difficult moving--boundary simulations, especially when the evolving geometry is strongly non-convex or contains rapidly changing features.

\section{Conclusion}\label{sec:conclusion}
An ALE--MFS framework for moving boundary problems in evolving two-dimensional domains has been investigated, with particular emphasis on the interaction between domain evolution, discretization, and \emph{a posteriori} error assessment in a meshless setting. Through numerical experiments of increasing geometric complexity and with fixed circular source placement, the general $\mathrm{LOOCV}\!-\!\ell_\infty$ indicator, computed via the hat matrix, was found to be a practical and reliable tool for parameter selection and error estimation: in well-resolved square MFS systems it provided a robust upper bound on the boundary error, and its decay under refinement reflected the expected convergence behaviour.

For highly non-convex and irregular geometries (such as the amoeba domain), the limitations of the classical square formulation became apparent. The interior point distribution deteriorated over time, and the LOOCV indicator became overly pessimistic and erratic, losing much of its diagnostic value. This motivated a zero-padded square construction that preserves $N_c = N_s$ while augmenting the system with zero rows to realise a stabilised least-squares solve without changing the source or collocation sets. In this setting, the general LOOCV indicator remained effective, conservative, and well correlated with the maximum-principle boundary error, which served as a complementary diagnostic based on the harmonic maximum principle.

A pseudoinverse-based pinv--Rippa algorithm was also tested as a low-cost extension of the classical Rippa formula. Because the original identity does not extend to zero-padded least-squares systems, this heuristic lacks theoretical guarantees and, in the numerical experiments, displayed irregular behaviour, especially on complex geometries. Consequently, it is best viewed as an auxiliary check, with $\mathrm{LOOCV}\!-\!\ell_\infty$ (in its hat-matrix form) and the maximum-principle indicator providing the primary \emph{a posteriori} assessment.

Beyond boundary-error analysis, the B-spline based ALE--MFS mesh motion was embedded into a finite-element setting and compared directly with a classical ALE--FEM mesh motion driven by the Barrett--Garcke--N\"urnberg curvature-flow scheme. Starting from identical initial triangulations and enforcing the same curvature-driven normal velocity, standard mesh-quality indicators were monitored: the minimum interior angle and the mesh ratio. For both the sharp asterisk and the highly non-convex amoeba geometries, and for coarse and fine meshes, the ALE--MFS motion systematically maintained larger minimum angles and smaller mesh-ratio growth than the ALE--FEM variant. In particular, for the amoeba geometry the ALE--FEM meshes approached typical remeshing thresholds, whereas the ALE--MFS meshes stayed further from these critical values. Since the ALE--MFS harmonic extension can be used as a drop-in replacement for the finite-element harmonic extension in existing ALE codes, this suggests a practical route to improving mesh robustness and delaying or avoiding costly remeshing steps.

Taken together, these findings indicate that ALE--MFS, equipped with solid \emph{a posteriori} indicators (in particular the general $\mathrm{LOOCV}\!-\!\ell_\infty$ and the maximum-principle error) and supported by inexpensive heuristics when desired, is a reliable and adaptable method for tackling moving boundary problems in challenging geometries. The combination of meshless harmonic extension, B-spline based curvature and normal reconstruction on point clouds, and rigorous error diagnostics makes the framework particularly attractive for strongly non-convex interfaces and large deformations, where mesh distortion is a primary concern.

Several directions for further study are natural. On the analytical side, a more detailed spectral analysis of MFS matrices under curvature-driven boundary data would clarify when square collocation is adequate and when least-squares regularisation becomes essential. On the algorithmic side, extending the approach to three dimensions, exploring geometry-adapted source placement strategies (such as normal-offset or multi-layer source sets), and combining ALE--MFS with higher-order or adaptive time-stepping schemes are promising avenues. Finally, applying the framework to fully coupled multiphysics problems, including free-surface Navier--Stokes and fluid--structure interaction, would provide stringent tests of its robustness and efficiency in realistic moving-interface simulations.

\section*{Acknowledgment}
This work was supported by the General Research Fund (GRF No. 12301824, 12300922) of the Hong Kong Research Grant Council and the Guangdong and Hong Kong Universities “1+1+1” Joint Research Collaboration Scheme  (project no. 2025A0505000014).

  \bibliographystyle{elsarticle-num}
  \bibliography{ref-abrv}
\end{document}